\documentclass {amsart}

\usepackage{setspace}

\usepackage{amsmath}
\usepackage{amssymb}
\usepackage{amsthm}
\usepackage{graphicx}

\input xy
\xyoption{all}

\newcommand{\mc} [1]{\mathcal{#1}}

\theoremstyle{plain}
\newtheorem{theorem}{Theorem}[section]
\newtheorem{lemma}[theorem]{Lemma}
\newtheorem{ellemma}[theorem]{Elementary Lemma}
\newtheorem{cor}[theorem]{Corollary}
\newtheorem{prop}[theorem]{Proposition}
\newtheorem{conjecture}[theorem]{Conjecture}

\theoremstyle{definition}
\newtheorem{definition}[theorem]{Definition}
\newtheorem{example}[theorem]{Example}

\newtheorem{remark}[theorem]{Remark}

\begin{document}

\title{A step in Castelnuovo theory via Gr\"obner bases}
\author {Ivan Petrakiev}
\address{Department of Mathematics, Harvard University, Cambridge, MA 02138}
\email{petrak@math.harvard.edu}
\begin {abstract} We establish the first previously unknown case of the Eisenbud-Harris conjecture in Castelnuovo theory concerning algebraic curves of high genus in ${\bf P}^n$. The problem is reduced to a question about zero-dimensional schemes $\Gamma \subset {\bf P}^{n-1}$ in symmetric position with certain constrains on the Hilbert function. The method of Gr\"obner bases is then applied to study the homogeneous ideal of $\Gamma$.
\end{abstract}
\thanks{The author was partially supported by the NSF Graduate Research Fellowship.}
\maketitle

\tableofcontents

\section{Introduction}

Let $C$ be a reduced, irreducible and nondegenerate curve of degree $d$ and arithmetic genus $g$ in ${\bf P}^{n}$, $n \geq 2$. A celebrated theorem of Castelnuovo (1889) gives an explicit upper bound $\pi_0 (d,n)$ on $g$ in terms of $d$ and $n$. Moreover, curves that attain the maximal genus, the so called {\it Castelnuovo curves}, have rather special extrinsic properties and are well understood. In particular, as long as $d \geq 2n+1$, such curves always lie on surfaces of minimal degree $n-1$. 

Castelnuovo's theorem has been reconsidered and extended further by several classical geometers, including G. Halph\'en (\cite{Halphen}), G. Fano (\cite{Fano}) and, much later, by Eisenbud-Harris (\cite{EH}). The main philosophy of the modern Castelnuovo theory is that curves of sufficiently high genus should lie on surfaces (or higher dimensional varieties) of some small degree. 

Extending Castelnuovo's bound, Eisenbud-Harris (\cite{EH}) defined a decreasing string of numbers $$\pi_\alpha(d,n) \approx \frac {d^2} {2(n-1+\alpha)} + O(d),$$ where $\alpha = 0,1,\dots, n-1,$ and made a conjecture: if $C$ is a curve of genus $g > \pi_\alpha(d,n)$ and $d \geq 2n+2\alpha-1$, then $C$ must lie on a surface of degree at most $n+\alpha-2$ (see Conjecture \ref{conj_1}).  

In \cite{EH}, a proof is given for the case $\alpha = 1$, although a similar result has been already known to Fano. The Eisenbud-Harris conjecture is also known to be true any $\alpha$, as long as $d >> 0$ (the explicit bound on $d$ is exponential in $n$).
\bigskip

In this paper we settle the next case $\alpha = 2$ ($n \geq 8$) of the Eisenbud-Harris conjecture (see thm. \ref{thm_cast}). The only previous work in this direction known to us is the paper of C. Ciliberto (\cite{Ci}) 
, where some partial results were obtained by different methods. 
\bigskip

Recall the main circle of ideas involved in Castelnuovo theory. Let $\Gamma = C \cap {\bf P}^{n-1}$ be a general hyperplane section of $C$. We will say, that $\Gamma$ is in {\it symmetric position} (see Def. \ref{def_sym}), which generalizes the notion of {\it uniform position}, first introduced by Harris in \cite{H}. As Castelnuovo observed, if $C$ is to have high genus, then $\Gamma$ must have a ``small'' Hilbert function $h_\Gamma(l)$ and, in particular, $\Gamma$ must fail to impose many conditions on quadrics in ${\bf P}^{n-1}$. Assume $d \geq 2n+1$. Then, according to the well-known Castelnuovo's lemma, $h_\Gamma(2)$ takes its minimal value $2n-1$ precisely when $\Gamma$ is a set of points lying on a {\it rational normal curve} in ${\bf P}^{n-1}$. This allowed Castelnuovo to determine his bound $\pi_0(d,n)$ on the genus of $C$ and describe the curves that achieve it (see Section \ref{sec_castelnuovo}). 

By generalizing Castelnuovo's lemma, one is naturally lead to conjecture the following: if $\Gamma \subset {\bf P}^{n-1}$ is a set of $d \geq 2n+2m-1$ points in symmetric position, with $h_\Gamma(2) \geq 2n+m-2$ (where $1\leq m \leq n-3$), then $\Gamma$ must lie on a curve of degree at most $n+m-2$ (see Conjecture \ref{conj_2}). 
\bigskip

In this paper we establish the first previously unknown cases $m=3$ ($n\geq 5$) and $m=4$ ($n\geq 7$) of Conjecture \ref{conj_2} (see thm. \ref{thm_cast2}). This result in turn implies thm. \ref{thm_cast}.
\bigskip

The starting point in our work is the fact, that under the assumptions of Conjecture \ref{conj_2}, $\Gamma$ lies on an {\it $m$-fold rational normal scroll} (the construction, due to Fano and, independently, Eisenbud-Harris, is described in Section \ref{sec_cast}). We use this, together with the symmetry of $\Gamma$, to write the beginning of a Gr\"obner basis for the homogeneous ideal of $\Gamma$ in degree 2, in a suitable coordinate system and monomial order (see Section \ref{sec_borel}-\ref{sec_grob}). It turns out, that there are only few quadrics missing in our Gr\"obner basis, precisely $\binom {m-1} 2$. We make a conjecture about the ``missing'' $\binom{m-1} 2$ quadrics and support it with some evidence, that comes from an elementary geometric observation (see lemma \ref{lemma_tfae}, lemma \ref{lemma_mult} and cor. \ref{cor_evidence}). In the cases $m=3,4$, we are actually able to complete the whole Gr\"obner basis of $\Gamma$ in degree 2, which allows us to prove our main results (in Section \ref{sec_points2}).

\subsubsection*{Notation and conventions.}

We work over algebraically closed field $\overline{k}$ of characteristic $0$. For any closed subscheme $X \subset {\bf P}^n$, denote by $$h_X(l) = \text{rk} \left( \rho: H^0 ({\bf P}^n, \mc{O}_{{\bf P}^n}(l)) \rightarrow H^0 (X,\mc{O}_X(l) \right)$$
the {\it Hilbert function} of $X$, where $\rho$ is the natural restriction map. Equivalently, $h_X(l)$ is the number of conditions that $X$ imposes on hypersurfaces of degree $l$ in ${\bf P}^n$.
We say, that $X \subset {\bf P}^n$ is {\it nondegenerate} if and only if $X$ is not contained in any hyperplane.
For any line bundle $\mc{L}$ on $X$ and linear system $V \subseteq |\mc{L}|$, denote by $Bs(V)$ the scheme-theoretic base locus of $V$.
If $X$ is integral, denote by $F(X)$ the field of functions on $X$.

\subsubsection*{Acknowledgements}
First, I would like to express my deepest gratitude to my thesis advisor Joe Harris for suggesting the problem and for guiding my research through the last 3 years. I would also like to thank David Eisenbud for a discussion that contributed to this paper.

\section {Projective curves with high genus}

\subsection {Castelnuovo's classical theorem} \label{sec_castelnuovo} \

\bigskip

A theorem of Castelnuovo provides an explicit bound on the genus $g$ of a nondegenerate curve of fixed degree $d$ in ${\bf P}^n$. Moreover, curves that attain the maximal genus, the so called {\it Castelnuovo curves}, always lie on surfaces of minimal degree (assuming $d \geq 2n+1$) and, therefore, are well understood. In this section we overview the basic ideas leading to Castelnuovo's classical result (thm. \ref{thm_cast1} below). These ideas will be extended further in the next section, where we discuss the refinements of Castelnuovo theory, due to Eisenbud-Harris. 

Let $C$ be a reduced, irreducible and nondegenerate curve in ${\bf P}^n$ of degree $d$ and genus $g$. Let $\Gamma = C \cap H$ be any hyperplane section of $C$. Then, it is a standard fact (\cite{EH}, Lemma 3.1), that
\begin{align}\label{eq_hyper}
h_{\Gamma}(l) \leq h_C (l) - h_C (l-1).
\end{align}

By the Riemann-Roch theorem, $h_C(l) = k d - g + 1,$ for any $k$ sufficiently large. By summing (\ref{eq_hyper}) over all $k \geq 1$, we get (\cite{EH}, Cor. 3.2):
\begin{align}\label{ineq_g}
g \leq \sum_{l=1}^\infty (d - h_\Gamma(l)).
\end{align}
Equality is achieved if and only if $C$ is Arithmetically Cohen-Macaulay (ACM) (\cite {Eis}). 

Inequality (\ref{ineq_g}) is central to Castelnuovo theory, since it gives an upper bound on the genus of the curve $C$ in terms of the Hilbert function of its hyperplane section $\Gamma$. In particular, $C$ has a maximal genus if and only if $C$ is ACM and $\Gamma$ has a minimal possible Hilbert function. 

We introduce the following two notions.

\begin{definition}\label{def_sym} Let $\Gamma_K \subset {\bf P}^{n-1}_K$ be an integral, nondegenerate zero-dimensional scheme of degree $d$, defined over a (non-algebraically closed) field $K$. Suppose, that $F(\Gamma_K) / K$ is a Galois extension with Galois group $S_d$, the full symmetric group on $d$ letters. Then, for any field extension $L/K$, such that the pull-back $\Gamma_L = \Gamma_K \times_K L \subset {\bf P}^{n-1}_L$ splits as a set of $d$ distinct (geometric) points, we say that $\Gamma_L$ is a set of points in {\it symmetric position} in ${\bf P}^{n-1}_L$.
\end{definition}

\begin{definition} (\cite{H}) Let $\Gamma \subset {\bf P}^{n-1}$ be a nondegenerate set of $d$ distinct points. Suppose, that any two subsets $A, B \subset \Gamma$ of the same cardinality also have the same Hilbert function, i.e. $h_A = h_B$. Equivalently, for any $l$, any hypersurface of degree $l$ containing $\geq h_\Gamma(l)$ points of $\Gamma$ necessarily contains all of $\Gamma$. Then, we say that $\Gamma$ is a set of points {\it in uniform position} in ${\bf P}^{n-1}$.
\end{definition}
Clearly, the following implications are valid: 
$$\text {(Symmetric position)} => \text {(Uniform position)} => \text{(Linearly general position)}$$

Let $C \subset {{\bf P}^n}_{\overline k}$ be a curve as above and let $H_K \cong {\bf P}^{n-1}_K$ be the generic hyperplane, defined over a pure transcendental field extension $K/\overline{k}$ of degree $n$. Let $\Gamma_K = C \cap H_K$ and let $L/K$ be any field extension, for which $\Gamma_L$ splits as a set of $d$ distinct points. By a standard monodromy argument (\cite{ACGH}), $\Gamma_L \subset H_L$ is a set of $d$ points in symmetric (hence uniform) position. 

In future, we will omit the subscript $L$ and will simply write $\Gamma = C \cap H \subset H$.

The Hilbert function of a set of points in uniform position satisfies the following ``subadditivity'' property (\cite{EH}, p.86):
\begin{align}\label{sub_add}
h_\Gamma(k+l) \geq \min \{ d, h_\Gamma(k) + h_\Gamma(l) - 1\}.
\end{align}

Now, since $\Gamma\subset {\bf P}^{n-1}$ is nondegenerate, we have $h_\Gamma(1) = n$. By (\ref{sub_add}), we have $h_\Gamma(2) \geq \min \{d, 2n-1 \}$. The following well-known lemma describes $\Gamma$ in the case when $h_\Gamma(2)$ is minimal, assuming that $d$ is not too small:

\begin{lemma}\label{lemma_castel}(Castelnuovo \cite{C}) Let $\Gamma$ be a set of $d \geq 2n+1$ points in linearly general position in ${\bf P}^{n-1}$. Then $h_\Gamma(2) \geq 2n-1$. Moreover, $h_\Gamma(2) = 2n-1$ if and only if $\Gamma$ lies on a rational normal curve.
\end{lemma}

We will present a proof of a generalized version of Castelnuovo's lemma later in Section \ref{sec_cast}.

\begin{remark}\label{bezout} Observe, that if $\Gamma$ lies on a rational normal curve $D$ then, by Bezout's theorem, any quadric in ${\bf P}^{n-1}$ containing $\Gamma$ also contains $D$. In particular, $h_\Gamma(2) = 2n-1$, which is the ``if'' part of Castelnuovo's lemma.
\end{remark}

Finally, we are in position to give the Castelnuovo's bound on the genus of space curves. For $d$, $g$ and $n$ as above, set $$\lambda_0 = \left[ \frac{d-1}{n-1}\right] \ \ \ \text{and} \ \ \ \epsilon_0 = (d-1) - \lambda_0 (n-1).$$
Define
\begin{align}\label{h_zero}
h_0(l) = 
\left\{
\begin{array}{ll}
l(n-1) +1, & 1 \leq l \leq \lambda_0 \\
d, & l \geq \lambda_0 +1 \\
\end{array}
\right.
\end{align}
and 
\begin{align}\label{pi_zero}
\pi_0 (d,n) = \sum_{l=1}^{\infty} (d - h_0 (l)).
\end{align}
Explicitly,
\begin{align}
\pi_0 (d,n) = \binom {\lambda_0} 2 (n-1) + \lambda_0 \epsilon_0.
\end{align}
We have (\cite{EH}, p.87):
\begin{theorem} (Castelnuovo)\label{thm_cast1} Let $C$ be a nondegenerate, reduced and irreducible curve of degree $d$, genus $g$ in ${\bf P}^n$. Then, $g \leq \pi_0 (d,n)$. Moreover, if $g = \pi_0 (d,n)$ and $d \geq 2n+1$, then $C$ lies on a surface $S$ of minimal degree $n-1$.
\end{theorem}

{\it Sketch of proof.}
Let $\Gamma = C \cap H$ be a general hyperplane section. From (\ref{sub_add}) and the fact that $h_\Gamma(1) = n$, we conclude that $h_\Gamma (l) \geq h_0 (l)$ for any $l$. By (\ref{ineq_g}) and (\ref{pi_zero}), we have $g \leq \pi_0 (d,n).$

Suppose that $d \geq 2n+1$ and $g = \pi_0 (d,n)$. It follows, that $h_\Gamma(l) = h_0(l)$ for any $l$. In particular, $h_\Gamma(2) = h_0(2) = 2n-1$. By Castelnuovo's lemma, $\Gamma$ lies on a rational normal curve $D \subset H$. It follows, that $C$ lies on a surface $S$ of minimal degree $n-1$ with $S \cap H = D$ (see \cite{EH}, p.87).


\subsection {Refinements due to Eisenbud-Harris}\label{sec_ref}\

\bigskip

It was recognized by Halphen, Fano and later by Eisenbud-Harris, that Castelnuovo's approach can be applied to study curves of fixed degree $d$ and nearly maximal genus $g$. The main hope is that such curves should lie on surfaces of some small degree (we will make this more precise later). In this section, we overview some results and conjectures in Castelnuovo theory, following Eisenbud-Harris  (\cite{EH}).


By taking a general hyperplane section of a curve $C$ with ``high'' genus, one is naturally lead to studying sets of points $\Gamma \subset {\bf P}^{n-1}$ in symmetric (or uniform) position with  ``small'' Hilbert functions. To make the problem more tractable, one may use the subadditivity inequality to bound the entire Hilbert function of $\Gamma$ just by knowing  $h_\Gamma(2)$ (of course, this bound need not be sharp in general). Again, the main hope is that if $h_\Gamma(2)$ is sufficiently small, that is, if $\Gamma$ fails to impose many conditions on quadrics in ${\bf P}^{n-1}$, then this failure is caused by the fact that $\Gamma$ lies on a curve $D$ of some small degree (compare to rmk. \ref{bezout}). 

In \cite{EH}, Eisenbud-Harris proved the following analogue of Castelnuovo's lemma:

\begin{lemma}\label{lemma_eh} (Eisenbud-Harris\footnote {According to \cite{C}, this lemma was already known to G. Fano.}) Let $n \geq 4$ and let $\Gamma$ be a set of $d \geq 2n+3$ points in uniform position in ${\bf P}^{n-1}$. If $h_\Gamma(2) = 2n$, then $\Gamma$ lies on an elliptic normal curve $D$ in ${\bf P}^{n-1}$.  
\end{lemma}

Also, Eisenbud-Harris introduced a new constant $\pi_1 (d,n) < \pi_0 (d,n)$ (see below) and, using the lemma above, proved:

\begin{theorem}\label{thm_eh} (\cite{EH}, p. 99) Let $n\geq 4$ and let $C \subset {\bf P}^n$ be a reduced, irreducible and nondegenerate curve of degree $d$ and genus $g$. Then,
\begin{enumerate}
\item If $g > \pi_1 (d,n)$ and $d \geq 2n+1$, then $C$ lies on a surface $S$ of degree $n-1$.

\item If $g = \pi_1 (d,n)$ and $d \geq 2n+3$, then $C$ lies on a surface $S$ of degree $\leq n$.
\end{enumerate}
\end{theorem}


They went even further to define a decreasing string of numbers $\pi_\alpha(d,n)$, where $0 \leq \alpha \leq n-1$ (see below), and proved:

\begin{theorem}(Eisenbud-Harris)\label {conj_eh_thm} For any $n$, there is a constant $d_0 = d_0(n)$ with the following property. Let $C$ be a reduced, irreducible, nondegenerate curve of genus $g$, degree $d \geq d_0$ in ${\bf P}^n$. If $g > \pi_\alpha (d,n)$, then $C$ lies on a surface of degree at most $n+\alpha-2$.
\end{theorem}

Eisenbud-Harris gave an explicit value of $d_0 (n)$ for which the theorem is known to be true. Unfortunately, their $d_0 (n)$ has exponential growth in $n$ (for example, $d_0(n) = 2^{n+1}$ for $n \geq 8$) and, obviously, is far from being sharp. In fact, they made the following:

\begin{conjecture}(Eisenbud-Harris)\label {conj_1} Let $C$ be a reduced, irreducible, nondegenerate curve of genus $g$, degree $d \geq 2n+2\alpha-1$ in ${\bf P}^n$. If $g > \pi_\alpha (d,n)$, then $C$ lies on a surface of degree at most $n+\alpha-2$.
\end{conjecture}

For example, the case $\alpha = 0$ is Castelnuovo's bound, while the case $\alpha = 1$ is part (a) of Thm. \ref {thm_eh}. 

One of the goals of this paper is to establish the next case $\alpha = 2$ ($n\geq 8$) of the Eisenbud-Harris conjecture, and give new partial results for $\alpha = 3,4$ (See Theorem \ref{thm_cast}). 

Conjecture \ref{conj_1} is closely related to the following conjecture for sets of points in symmetric position with small $h_\Gamma(2)$:

\begin{conjecture}\label {conj_2} Let $\Gamma$ is a set of $d \geq 2n+2m-1$ points in symmetric position in ${\bf P}^{n-1}$, where $1 \leq m \leq n-3$. Suppose, that $h_\Gamma(2) = 2n+m-2$. Then, $\Gamma$ lies on a curve $D$ of degree at most $n+m-2$. 
\end{conjecture}
For example, the case $m=1$ is Castelnuovo's lemma, while the case $m=2$ is the Eisenbud-Harris lemma. 

One of the main goal of this paper is to establish the cases $m=3$ ($n\geq 6$) and $m=4$ ($n \geq 8$) of the conjecture above (see Theorem \ref{thm_cast2}).

\begin{remark}\label{rmk_comp} As we will see in Section \ref{sec_grob}, if the curve $D$ exists, then $D$ is unique and $D$ is a component of $Bs(|\mc{I}_\Gamma(2)|)$.
\end{remark}

\begin{remark} The assumption $m \leq n-3$ is necessary. For example, if $m=n-2$, then one could take $\Gamma$ to be a complete intersection of a del Pezzo surface (of degree $n$) and two general quadrics to produce a counterexample.
\end{remark}

\begin{remark} Conjecture \ref{conj_2} is a very special case of a conjecture in so called Higher Castelnuovo Theory, introduced by Eisenbud-Green-Harris in \cite{EGH}. In their paper, the conjecture is stated with the weaker assumption that $\Gamma$ is in {\it uniform} position. In our work, the assumption that $\Gamma$ is in symmetric position will play an essential role.
\end{remark}

\subsection{The definition of $\pi_\alpha(d,n)$}\label{sec_defpi}\

\bigskip

We recall some definitions from \cite{EH}.

Let $C$ be a reduced, irreducible curve of degree $d$ and genus $g$ in ${\bf P}^n$. For any $0 \leq \alpha \leq n-2$, set $$\lambda_\alpha = \left[ \frac {d-1} {n-1+\alpha} \right] \ \ {\text {and}} \ \ \epsilon_\alpha = d-1-\lambda_\alpha(n-1+\alpha).$$
Also, set
$$\mu_\alpha = \max \left(0, \left[ \frac {\alpha-n+2+\epsilon_\alpha} 2 \right] \right).$$
Define the function $h_\alpha(l)$, depending also on $d$ and $n$, as follows:
\begin{align}
h_\alpha(l) = 
\left\{
\begin{array}{ll}
l(n+\alpha-1) - \alpha+1, & 1 \leq l \leq \lambda_\alpha \\
d - \mu_\alpha, & l = \lambda_\alpha +1 \\
d, &  l \geq \lambda_\alpha + 2 \\
\end{array}
\right.
\end{align}
Finally, define
\begin{align}\label{eqn_pi}
\pi_\alpha(d,n) = \sum_{l=1}^\infty (d- h_\alpha(l)).
\end{align}
Explicitly,
\begin{align}
\pi_\alpha(d,n) = \binom{\lambda_\alpha} 2 (n-1+\alpha) + \lambda_\alpha(\epsilon_\alpha+\alpha) + \mu_\alpha.
\end{align}
Asymptotically, we have $\pi_\alpha(d,n) \approx \frac {d^2} {2(n-1+\alpha)} + O(d).$ 

Similar definitions exist in the case $\alpha = n-1$ (see \cite{EH}). 

The relevance of the definitions above comes from the following lemma, which is a slight refinement of \cite{EH}, Lemma 3.24(i).

\begin{lemma}\label{lemma_normal} 
Let $\Gamma \subset {\bf P}^{n-1}$ be a general hyperplane section of $C$. Suppose, that $\Gamma$ lies on an irreducible curve $D \subset {\bf P}^{n-1}$ of degree $n+\alpha-1$, for some $0 \leq \alpha \leq n-2$. Assume, that $D$ is $(\lambda_\alpha+1)$-normal, that is, the linear system on $D$ cut out by hypersurfaces of degree $\lambda_\alpha+1$ is complete. Then, $h_\Gamma(l) \geq h_\alpha(l)$ for all $l$.
\end{lemma}

The proof is the same as in \cite{EH}, so we omit it.

\begin{remark}\label{rmk_gr} By the result of Gruson, et.al. \cite{Gr}, $D$ is automatically $(\lambda_\alpha+1)$-normal in the following cases: (i) $\lambda_\alpha \geq \alpha$, or (ii) $\lambda_\alpha = \alpha-1$ and $D$ is not a smooth rational curve with an $(\alpha+2)$-secant line.  
\end{remark}


\subsection {New results for $\alpha\leq 4$}\label{new_results}\

\bigskip

Our new results towards Conjecture \ref{conj_1} are the following.

\begin{theorem}\label{thm_cast} Let $n\geq 8$ and let $C$ be a reduced, irreducible and nondegenerate curve of degree $d$ and genus $g$ in ${\bf P}^n$. Then:
\begin{enumerate}
\item If $g > \pi_2 (d,n)$ and $d \geq 2n+3$, then $C$ lies on a surface $S$ of degree $\leq n$.
\item If $g > \pi_3 (d,n)$ and $2n+5 \leq d \leq 4n+6$, then $C$ lies on a surface $S$ of degree $\leq n+1$.
\item If $g > \pi_4(d,n)$ and $2n+7 \leq d \leq 3n+3$, then $C$ lies on a surface $S$ of degree $\leq n+2$.
\end{enumerate}
\end{theorem}

The proof relies on knowing Conjecture \ref{conj_2} for $m=3,4$. We have:
\begin{theorem}\label{thm_cast2} Let $\Gamma \subset {\bf P}^{n-1}$ be a set of $d $ points in symmetric position. 
\begin{enumerate}
\item Let $n \geq 5$. If $h_\Gamma(2) = 2n+1$ and $d \geq 2n+5$, then $\Gamma$ lies on a curve $D$ of degree $\leq n+1$.
\item Let $n \geq 7$. If $h_\Gamma(2) = 2n+2$ and $d \geq 2n+7$, then $\Gamma$ lies on a curve $D$ of degree $\leq n+2$.
\end{enumerate}
\end{theorem}
The proof of Theorem \ref{thm_cast2} will occupy Sections \ref{sec_points1} and \ref{sec_points2} of this paper. Here, we merely show how one result implies the other. 
\bigskip

{\it Proof of Theorem \ref{thm_cast} assuming Theorem \ref{thm_cast2}:}

Let $g > \pi_\alpha(d,n)$ be as in the theorem, $\alpha  \leq 4$. We want to show, that $C$ lies on a surface $S$ of degree $\leq n+\alpha-2$. 

First, we reduce to the case when $C$ is linearly normal. If this is not the case, then $C$ is a projection of some curve $C' \subset {\bf P}^{n+1}$ of degree $d$ and genus $g$. Notice, that $\pi_{\alpha}(d,n) \geq \pi_{\alpha-1} (d,n+1)$. By induction on $\alpha$, we conclude that $C'$ lies on a surface $S'$ of degree $\leq n+\alpha-2$. We may now take $S$ to be the projection of $S'$ in ${\bf P}^n$.

So, assume that $C$ is linearly normal and let $\Gamma = C \cap H$ be a general hyperplane section. Our strategy is the following: use the assumption on the genus $g > \pi_\alpha(d,n)$ and the subadditivity of $h_\Gamma(l)$ to bound $h_\Gamma(2)$. Then, apply Theorem \ref{thm_cast2} to conclude, that $C$ lies on a curve $D$ of some ``small'' degree. By lemma \ref{lemma_normal}, conclude that, in fact, $\deg D \leq n+\alpha-2$. 

Once we know this, it follows at once that $C$ lies on a surface $S$ of degree $\leq n+\alpha-2$. Indeed, by linear normality of $C$, every quadric in $H$ containing $\Gamma$ is a restriction of a quadric in ${\bf P}^n$, containing $C$. Now, $D$ is a component of the base locus $Bs (|\mc{I}_\Gamma(2)|)$ in $H$ (see remark \ref{rmk_comp}). It follows, that $Bs(|\mc{I}_C(2)|)$ has a two-dimensional component $S$, whose general hyperplane section is $D$.  

a) Assume first, that $d \geq 2n+7$. From the assumption $g > \pi_2 (d,n)$ and the subadditivity of $h_\Gamma(l)$, it follows that $h_\Gamma(2) \leq 2n+2$.  By Theorem \ref{thm_cast2}, $\Gamma$ lies on a curve $D$ of degree $\leq n+2$. 

We claim, that $\deg D \leq n$. Indeed, the condition $g > \pi_2(d,n)$, together with (\ref{ineq_g}) and (\ref{eqn_pi}) implies, that
$h_\Gamma(l) \leq h_2 (l)$,
for some $l \geq 1$. By Lemma \ref{lemma_normal}, $D$ cannot be of degree $n+1$ (notice, that $\lambda_2(d,n) \geq 2$, so remark \ref{rmk_gr} applies). Similarly, $g > \pi_2 (d,n) \geq \pi_3 (d,n)$ implies that $h_\Gamma(l) \leq h_3(l)$ for some $l \geq 1$. By Lemma \ref{lemma_normal}, $D$ cannot be of degree $n+2$ (we have $\lambda_3 (d,n) \geq 2$ and $D$ is not smooth rational, because $h_D(2) = h_\Gamma(2) \leq 2n+2$; so, remark \ref{rmk_gr} still applies). We conclude, that $\deg D \leq n$.

Hence, $C$ lies on a surface $S$ of degree $\leq n$.

It remains to consider the case $2n+3 \leq d \leq 2n+6$. Now, the assumption on the genus $g$ implies the stronger condition $h_\Gamma(2) \leq 2n$. By the Eisenbud-Harris lemma, $\Gamma$ lies on a curve $D$ of degree $\leq n$, and we are in the same situation as above.

b, c) In each case, the bounds on $d$ are so chosen as to imply $h_\Gamma(2) \leq 2n+2$. The rest of the argument is very similar to part a) and is omitted.
\hfill$\Box$

\section {Points in symmetric position with small Hilbert functions} \label{sec_points1}

The following two sections will be devoted to proving our new results towards Conjecture \ref{conj_2} that were announced in Theorem {\ref{thm_cast2}} above.

Since we will work exclusively with sets of points in symmetric position, it will be notationally convenient to replace $n-1$ with $n$ everywhere. In particular, Conjecture \ref{conj_2} takes the form: 

\bigskip

{\bf Conjecture.} {\it Let $\Gamma$ is a set of $d \geq 2n+2m+1$ points in symmetric position in ${\bf P}^{n}$, where $1 \leq m \leq n-2$. Suppose, that $h_\Gamma(2) = 2n+m$. Then, $\Gamma$ lies on a curve $D$ of degree at most $n+m-1$. } 

\bigskip

In what follows, we start with some generalities on points in symmetric position, and then gradually specialize to the cases that are of interest to us. 

\subsection{Basic setup}\label{sec_borel}\

\bigskip

Let $\Gamma \subset {\bf P}^n$ be a set of $d$ points in {\it linearly general} position, defined over some field $L$, and let $I_\Gamma$ be the homogeneous ideal of $\Gamma$. Let $p_0,\dots,p_n$ be some of the points of $\Gamma$. Choose homogeneous coordinates $x_0,\dots,x_n$ in ${\bf P}^n$, for which $p_0,\dots,p_n$ are the standard vertices. 

We set {\it graded lexicographic order (grlex)} on the set of monomials in the polynomial ring $L[x_0,\dots,x_n]$ (see \cite{CLO} for definitions). For any $f \in L[x_0,\dots,x_n]$, denote by $in(f)$ the initial monomial of $f$.

Consider the homogeneous ideal $I_\Gamma$ and its initial ideal $in (I_\Gamma)$. Denote by $in(I_\Gamma)_l$ the $l$-th graded piece of $in(I_\Gamma)$.

Since $I_{p_i} = (x_0,...,\widehat{x_i},...,x_n)$ for $i=0,\dots,n$, we have $x_i^2 \notin in(I_\Gamma) \subset in(I_{p_i})$.

\begin{prop}\label{prop_borel} Assume that $\Gamma$ is in symmetric position. Then, $in(I_\Gamma)$ has the following property: if $x_i x_j \in in(I_\Gamma)$ (where $i<j$), then $x_{i-1} x_{j} \in in(I_\Gamma)$ and $x_i x_{j-1} \in in(I_\Gamma)$, unless $i=0$ or $i = j-1$, respectively.
\end{prop}

\begin{proof} Let $Q \in in(I_\Gamma)$ be a quadric with $in(Q) = x_i x_j$. Geometrically, this means that: (i) The singular locus $sing(Q)$ contains the linear span $\overline{p_0\dots p_{i-1}}$, but not $\overline{p_0\dots p_i}$; (ii) $Q$ contains the lines $\overline{p_0 p_{1}}, \dots,\overline{p_0 p_{j-1}}$, but not $\overline {p_0 p_{j}}$.

Suppose that $i > 0$ and consider the permutation $\tau' = \tau_{i-1,i}$ on $\{p_0,\dots,p_n\}$ that transposes $p_{i-1}$ and $p_i$. Since $\Gamma$ is in symmetric position, there is a quadric $Q' = \tau' (Q) \in I_\Gamma$, such that: (i) $Q'$ is singular at $\overline{p_0\dots p_{i-2} p_i}$, but not $\overline{p_0\dots p_{i-1} p_i}$; (ii) $Q'$ contains the lines $\overline{p_0 p_1},\dots,\overline{p_0 p_{j-1}}$, but not $\overline{p_0 p_j}$. It follows that $in(Q') = x_{i-1} x_j$.

Similarly, suppose that $i < j-1$ and consider the transposition $\tau'' = \tau_{j-1,j}$. By a similar argument, we have $in(Q'') = x_i x_{j-1}$, where $Q'' = \tau''(Q)$.
\end{proof}

From now on, $\Gamma \subset {\bf P}^n$ is always assumed to be in symmetric position.

For any $k = 1,\dots,n$, we may consider the projection
$$\pi_k : {\bf P}^n \dashrightarrow {\bf P}^k$$
from the linear subspace $\{x_{n-k} = \dots = x_{n} = 0\} \cong {\bf P}^{n-k-1}$. Let
$$\Gamma_k = \pi_k (\Gamma - \{p_0,\dots,p_{n-k-1}\}).$$
By the well-known {\it Elimination theorem} (\cite{CLO}), the initial ideal of $\Gamma_k$ is:
$$in(I_{\Gamma_k}) = in(I_\Gamma) \cap k[x_{n-k},\dots,x_n].$$
It will be useful to depict $in(I_\Gamma)_2$ in diagrams of the form:

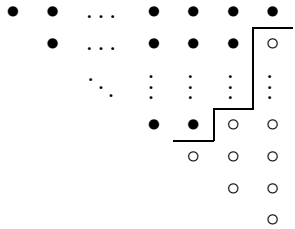
\begin{singlespace}
\begin{align*}
\begin{array}{r}
	\begin{array}{lllllll}
	\bullet & \bullet &\dots  & \bullet & \bullet & \bullet & \bullet \\
	\end{array} \\
	\begin{array}{llllll|l}
		\cline{7-7}
       		& \bullet &\dots  & \bullet & \bullet & \bullet  & \circ \\
        	&         &\ddots  & \vdots    & \vdots    & \vdots   & \vdots\\
	\end{array} \\
	\begin{array}{lllll|lll}
		\cline{6-6}
	        &         &       & \bullet  & \bullet & \circ & \circ \\
	\end{array} \\
	\begin{array}{lllllll}
		\cline{5-5}
        	&         &       &               & \circ & \circ   & \circ \\
        	&         &       &               &         & \circ   & \circ \\
        	&         &       &               &         &         & \circ \\
	\end{array}
\end{array}
\end{align*}
\end{singlespace}
\begin{center}
{\bf An example of $in(I_\Gamma)_2$}
\end{center}
In the diagram above, every dot corresponds a monomial $x_i x_j$ with $i < j$. The monomials are ordered lexicographically (from left to right, top to bottom). In particular, the dots in the $(i+1)$-st row correspond to monomials:
$$\begin{array}{cccc}
x_i x_{i+1} & x_i x_{i+2} & \dots & x_i x_n \\
\end{array}
$$
By prop. \ref{prop_borel}, there is a step-like line separating the monomials in $in(I_\Gamma)_2$ from the other degree two monomials.

By the Elimination theorem, the diagram of $in(I_{\Gamma_k})_2$ can be obtained from  the diagram of $in(I_\Gamma)_2$ by deleting the top $n-k$ rows.

\begin{example} If $in(I_\Gamma)_2$ is as in the diagram above, then $in(I_{\Gamma_4})_2$ is:
\begin{singlespace}
\begin{align*}
\begin{array}{cccc}
\bullet &\bullet &\circ  &\circ   \\
        &\circ   &\circ  &\circ   \\
        &        &\circ  & \circ  \\
        &        &       & \circ  \\
\end{array}
\end{align*}
\end{singlespace}
\end{example}

\bigskip

\subsection{The Generalized Castelnuovo Lemma}\label {sec_cast}\

\bigskip

The following basic lemma plays a key role in Castelnuovo theory.

\begin{lemma}\label {lemma_cast} (Fano, Eisenbud-Harris) Let $\Gamma$ be a set of $d \geq 2n+2m+1$ points in symmetric position in ${\bf P}^n$, where $m \leq n-2$. If $h_\Gamma(2) = 2n+m$, then there exists an $m$-fold rational normal scroll $\Sigma$ containing $\Gamma$.
\end{lemma}
In particular, if we set $m=1$, we recover Castelnuovo's classical lemma.

\begin{remark} In fact, Lemma \ref{lemma_cast} is true if $\Gamma$ is just assumed to be in uniform position (\cite{EH}), but we will not use this in the sequel. 
\end{remark}

Below we sketch the construction of $\Sigma$. The actual proof can be found in \cite{EH}. 
\bigskip

{\it Step 1.} Let $\Gamma =\{p_0,\dots, p_{d-1}\}$. Choose homogeneous coordinates $x_0, \dots, x_n$ in ${\bf P}^n$, for which $p_0,\dots, p_n$ are the standard vertices. Let $\Lambda = \{ x_{n-1} = x_n = 0\}  \cong {\bf P}^{n-2}$, i.e. $\Lambda = \overline {p_0,\dots,p_{n-2}}$. Then, $\Lambda$ imposes at most $\binom {n-1} 2$ conditions on the linear system $|\mc{I}_\Gamma(2)|$. It follows that $h^0 (\mc{I}_{\Gamma \cup \Lambda} (2)) \geq n-m$. So, we may take $n-m$ (general) linearly independent quadrics
$$Q_i = \left|
\begin{array}{ll}
L_i & x_{n-1} \\
M_i & x_n \\
\end{array}
\right|
$$
vanishing on $\Gamma$, where $0 \leq i \leq n-m-1$.
\bigskip

{\it Step 2.} It is easy to show, that the linear forms $L_i$ (resp. $M_i$) are linearly independent on $\Lambda$ (see \cite{EH}). By symmetry, it follows that $L_i$'s (resp. $M_i$'s) are linearly independent on $p_0, \dots, p_{n-m-1}$. So, after a linear change, we may assume that $L_i$ vanishes on $p_0,\dots, \widehat{p_i}, \dots, p_{n-m-1}$, for $0 \leq i \leq n-m-1$, and $L_i(p_i) = 1$. It follows immediately that $M_i$ vanishes on $p_0,\dots, \widehat{p_i},\dots, p_{n-m-1}$, and $M_i(p_i) \neq 0$.
\bigskip

{\it Step 3.} Define the determinantal locus

\begin{align}\label {eq_sigma} \Sigma := \text{rk} \left|
\begin{array}{lllll}
L_0 & L_1 & \dots & L_{n-m-1} & x_{n-1} \\
M_0 & M_1 & \dots & M_{n-m-1} & x_n \\
\end{array}
\right| \leq 1
\end{align}

Geometrically, $\Sigma$ is just the residual intersection of the quadrics $Q_i$ with respect to $\Lambda$. It is explained in \cite{EH}, that $\Sigma$ is in fact a rational normal $m$-fold scroll.
\bigskip

{\it Step 4.} It remains to show that $\Sigma$ contains $\Gamma$. A priori, it is clear that $\Sigma$ contains $\Gamma - (\Lambda\cap \Gamma) = \{p_{n-1}, p_n,\dots, p_{d-1}\}$. 

Consider the quadrics 
$$Q_{ij} = \left|
\begin{array}{ll}
L_i & L_j \\
M_i & M_j \\
\end{array}
\right|
$$
for $0 \leq i < j \leq n-m-1$. Now, $Q_{ij}$ vanishes on $p_0,\dots, \widehat {p_i}, \dots, \widehat {p_j}, \dots, p_{n-m-1}$, and also on $p_{n-1}, p_n, \dots, p_{d-1}$. Since $Q_{ij}$ contains at least $$n-m-2 + (d-n+1) \geq 2n+m = h_\Gamma(2)$$ points of $\Gamma$, it follows by symmetry that $Q_{ij}$ contains all points of $\Gamma$. This completes the last step.

\begin{lemma}\label {lemma_rule} Let $\Gamma$ and $\Sigma$ be as in the previous lemma.

(a) $\Sigma_{\text{sing}} \cap \Gamma = \{\emptyset \}$

(b) No two points of $\Gamma$ lie in the same ruling of $\Sigma$.
\end{lemma}

\begin{proof} Let $S_d$ be the group of permutations on $d$ letters $\{p_0,\dots,p_{d-1}\}$, and let $G$ be the subgroup of $S_d$ that leaves the two disjoint subset $\{p_0,\dots,p_{n-2}\}$ and $\{p_{n-1},\dots,p_{d-1}\}$ invariant. Then, the construction of $\Sigma$ in lemma \ref{lemma_cast} depends on the order of points $\Gamma = \{p_0,\dots,p_{d-1}\}$ only upto a $G$-action.

(a) Since $\Sigma \subset {\bf P}^n$ is an $m$-fold rational normal scroll, $m\leq n-2$, the singular locus $\Sigma_{sign}$ is a linear subspace of ${\bf P}^n$ of dimension $\leq n-4$. If a point $p \in \Gamma$ lies in the singular locus of $S$ then so does any point in the orbit $G \cdot p_i$. This contradicts the linear generality of $\Gamma$.

(b) A ruling of $\Sigma$ is a linear subspace of ${\bf P}^{n}$ of dimension $\leq n-3$. If two points $p,q \in \Gamma$ lie in the same ruling of $\Sigma$, then any point in the orbit $G \cdot p$ also lies in that ruling. Again, this contradicts the linear generality of $\Gamma$.
\end{proof}

Our next goal is to determine a Gr\"ober basis for $\Sigma$ for the lexicographic order on $x_0, \dots, x_n$. From {\it Step 2} in the construction of $\Sigma$, we have:
$$in (L_i) = in (M_i) = x_i$$
for $0 \leq i \leq n-m-1$. Therefore,
$$in (Q_i) = in (L_i x_n - M_i x_{n-1}) = x_{i} x_{n-1}.$$
By lemma \ref{lemma_rule}(b), any two points $p_i$, $p_j$ lie on different rulings of $\Sigma$. It follows that
$$in (Q_{ij}) = in (L_i M_j - L_j M_i) = x_i x_j.$$

\begin{lemma}\label {lemma_grob} The $\binom {n-m+1} 2$ quadrics $\{Q_i\} \cup \{ Q_{ij}\}$ defining $\Sigma$ form a Gr\"obner basis for $\Sigma$. In other words, $in(I_\Sigma)$ is generated by the monomials
$$\{ x_i x_j \}_{0 \leq i < j \leq n-m-1} \bigcup \{ x_i x_{n-1} \}_{0 \leq i \leq n-m-1}. \eqno(2)$$
\end{lemma}

\begin{proof} The rational normal scroll $\Sigma$ has Hilbert function
$h_\Sigma(t) = (n-m) \binom {t+m-1}{m} + \binom {t+m} m.$
It is easy to check that the ideal generated by the monomials in $(2)$ has the same Hilbert function. The lemma now follows from Macaulay's theorem.
\end{proof}

The initial ideal of $\Sigma$ in degree 2 is represented by the following diagram:
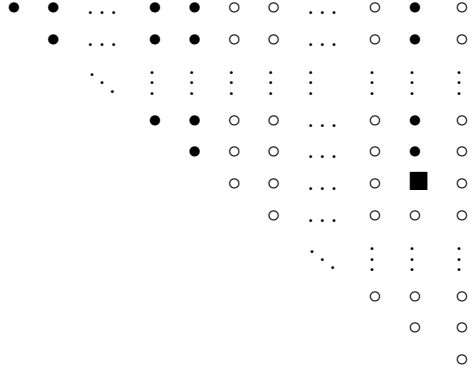
\begin{singlespace}
\begin{align*}
\begin{array}{lllllllllll}
\bullet & \bullet  & \dots & \bullet & \bullet & \circ & \circ & \dots & \circ & \bullet & \circ \\
        & \bullet  & \dots & \bullet & \bullet & \circ & \circ & \dots & \circ & \bullet & \circ \\
        &          & \ddots& \vdots  & \vdots  & \vdots& \vdots& \vdots& \vdots  & \vdots & \vdots \\
        &          &       & \bullet & \bullet& \circ  & \circ & \dots   & \circ & \bullet & \circ \\
        &          &       &         & \bullet & \circ & \circ & \dots & \circ & \bullet & \circ \\
        &          &       &         &         & \circ & \circ & \dots & \circ & \blacksquare&\circ \\
        &          &       &         &         &       & \circ & \dots & \circ & \circ   & \circ \\
        &          &       &         &         &       &       & \ddots& \vdots & \vdots  & \vdots \\
        &          &       &         &         &       &       &       & \circ & \circ   & \circ \\
        &          &       &         &         &       &       &       &       & \circ   & \circ \\
        &          &       &         &         &       &       &       &       &         & \circ \\
\end{array}
\end{align*}
\end{singlespace}
\begin{center}
{\bf The initial ideal of $\Sigma$}
\end{center}
\bigskip

The monomial $x_{n-m-1} x_{n-1} = in(Q_{n-m-1})$ plays a special role, and it is marked with $\blacksquare$ in the diagram above. 

\begin{lemma}\label {lemma_proj} Consider the projection $\pi_{m+1} : {\bf P}^n \dashrightarrow {\bf P}^{m+1}$. Then, the image $\pi(\Sigma)$ is the quadric $\pi (Q_{n-m-1})$. Moreover, the restricted map $\pi|_\Sigma : \Sigma \dashrightarrow \pi(\Sigma)$ is birational and the locus of indeterminancy of $\pi|_\Sigma$ consists precisely of the points $\{ p_0,\dots, p_{n-m-2} \}$.
\end{lemma}

\begin{proof} By elimination theory, we have:
$I_{\pi(\Sigma)} = I_\Sigma \cap k[x_{n-m-1},\dots,x_n].$ By lemma \ref{lemma_grob}, we see that $\pi(\Sigma) = \pi(Q_{n-m-1})$, which is a quadric in ${\bf P}^{m+1}$.
\end{proof}

\begin{cor}\label{cor_birat} For any $k \geq m+1$, let $\Sigma_k = \pi_k (\Sigma) \subset {\bf P}^k$. Then, $\Sigma_k \subset {\bf P}^k$ is a rational normal $m$-fold scroll, and the restricted map $\pi|_\Sigma : \Sigma \dashrightarrow \Sigma_k$ is birational.
\end{cor}

\subsection{On the initial ideal of $\Gamma$}\label{sec_grob}\

\bigskip

We adopt the setting from the previous section. In particular, let $\Gamma$ be a set of $d \geq 2n+2m+1$ points in symmetric position in ${\bf P}^n$ ($m \leq n-2$), such that $h_\Gamma(2) = 2n+m$. By lemma \ref{lemma_cast}, there is an $m$-fold rational normal scroll $\Sigma$ containing $\Gamma$. 

Since $\Gamma$ is in symmetric position and $I_\Gamma \supset I_\Sigma$, lemma \ref{lemma_grob} and prop. \ref{prop_borel} give:

\begin{cor}\label {cor_mono} Any monomial $x_i x_j$ with $0 \leq i \leq n-m-1$, $0 \leq j \leq n-1$ and $i < j$ occurs in the initial ideal of $\Gamma$.
\end{cor}

The monomials occurring in cor. \ref{cor_mono} are depicted in the following diagram:
\newpage
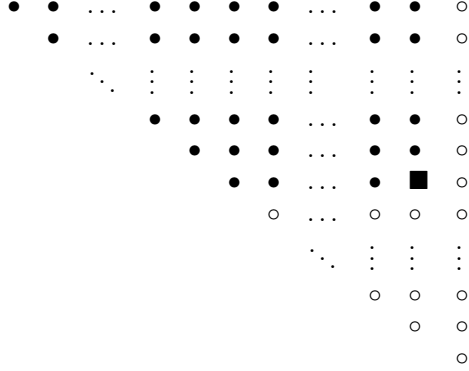
\begin{singlespace}
\begin{align*}
\begin{array}{lllllllllll}
\bullet & \bullet  & \dots & \bullet & \bullet & \bullet & \bullet & \dots & \bullet & \bullet & \circ \\
        & \bullet  & \dots & \bullet & \bullet & \bullet & \bullet & \dots & \bullet & \bullet & \circ \\
        &          & \ddots& \vdots  & \vdots  & \vdots & \vdots& \vdots& \vdots  & \vdots & \vdots \\
        &          &       & \bullet & \bullet& \bullet  & \bullet & \dots   & \bullet & \bullet & \circ \\
        &          &       &         & \bullet & \bullet & \bullet & \dots & \bullet & \bullet & \circ \\
        &          &       &         &         & \bullet & \bullet & \dots & \bullet & \blacksquare&\circ \\
        &          &       &         &         &       & \circ & \dots & \circ & \circ   & \circ \\
        &          &       &         &         &       &       & \ddots& \vdots & \vdots  & \vdots \\
        &          &       &         &         &       &       &       & \circ & \circ   & \circ \\
        &          &       &         &         &       &       &       &       & \circ   & \circ \\
        &          &       &         &         &       &       &       &       &         & \circ \\
\end{array}
\end{align*}
\end{singlespace}
\begin{center} {\bf Part of the initial ideal of $\Gamma$ in degree 2}
\end{center}
\bigskip

For any monomial $x_i x_j$ with $0 \leq i \leq n-m-1$ and $0 \leq j \leq n-1$, we can choose a representative $Q_{ij}$ in $H^0 (I_\Gamma(2))$ with $in(Q_{ij}) = x_i x_j$. Let $\mc{V} \subset H^0 (I_\Gamma(2))$ be the linear subspace spanned by the $Q_{ij}$'s. 

\begin{cor}\label{cor_codim} The codimension of $\mc{V}$ in $H^0 (\mc{I}_\Gamma(2))$ is $\binom {m-1} 2$. In particular, it does not depend on $n$.
\end{cor}

\begin{proof} We have $h^0 (\mc{I}_\Gamma(2)) = \binom {n+2} 2 - (2n+m)$ and $\dim \mc{V} = \binom{n-m} 2 + m(n-m)$.
Hence, $h^0 (\mc{I}_\Gamma(2)) - \dim \mc{V} = \binom {m-1} 2.$
\end{proof}

\begin{lemma}\label {prop_bp} For each $0 \leq i \leq n-m-1$, the Zariski tangent space $T_{p_i} Bs (|\mc{V}|)$ is 1-dimensional.
\end{lemma}

\begin{proof} By symmetry, it suffices to check the claim for $p_0$. Observe that the monomials $x_0 x_1, \dots, x_0 x_{n-1}$ occur as initial monomials in $\mc{V}$, while $x_0 x_n$ does not. The claim follows immediately.
\end{proof}

\begin{cor}\label{cor_hull} Suppose, that $\Gamma$ lies on a reduced, irreducible curve $D$ with $ \deg(D) < |\Gamma| / 2$. Then, $D$ is a component of $Bs(|\mc{I}_\Gamma(2)|)$. Moreover, $D$ is a unique curve with this property.
\end{cor}
\begin{proof} By Bezout's theorem, $D$ is contained in $Bs(|\mc{I}_\Gamma(2)|)$. By lemma \ref{prop_bp}, the Zariski tangent space of $Bs(|\mc{I}_\Gamma(2)|)$ at, say $p_0$, is at most 1-dimensional. The claim follows.
\end{proof}

Let us state the following weak version of Conjecture \ref{conj_2}:
\bigskip

{\bf Weak Conjecture \ref{conj_2}.} {\it For each $p_i \in \Gamma$, the Zariski tangent space $T_{p_i} Bs(|\mc{I}_\Gamma(2)|)$ is 1-dimensional.}
\bigskip
 
We have:

\begin{lemma}\label {lemma_tfae} The following are equivalent:
\begin{enumerate}
\item Weak Conjecture holds for $\Gamma$;
\item The monomial $x_0 x_n$ does not belong to $in (\Gamma)$;
\item None of the monomials $x_i x_n$ belongs to $in (\Gamma)$;
\item The following holds: $\dim \left( in(I_\Gamma)_2 \cap k[x_{n-m},\dots,x_n] \right) = \binom {m-1} 2$;
\item The following holds: $\dim \left( in(I_\Gamma)_2 \cap k[x_{n-m-1},\dots,x_n] \right) = \binom {m} 2 + 1$.
\end{enumerate}
\end{lemma}

\begin{remark}\label{remark_inq} A priori, we know that inequality $\leq$ holds in part (d) and (e) of the lemma above. This follows from cor. \ref{cor_codim}.
\end{remark}

\begin{proof} 

$(a) \Leftrightarrow (b)$ Same argument as in the proof of lemma \ref{prop_bp}.

$(b) \Leftrightarrow (c)$ follows from prop. \ref{prop_borel}.

$(c) \Leftrightarrow (d)$ follows from the description of the initial monomials in $\mc{V}$ (cor. \ref{cor_mono}) and from cor. \ref{cor_codim}.

$(d) \Leftrightarrow (e)$ follows from the description of the initial monomials in $\mc{V}$ and from prop. \ref{prop_borel}.
\end{proof}

Consider the projection $\pi_{m+1} : {\bf P}^n \dashrightarrow {\bf P}^{m+1}$ as in section \ref{sec_borel}, and let $\Gamma_{m+1}$ be the image of $\Gamma - \{p_0,\dots,p_{n-m-2}\}$. Let $\Gamma_{m+1} = \{q_0,\dots,q_{e-1}\}$, where $e = d - (n-m-1)$, with $q_i = \pi_{m+1} (p_{i+n-m-1})$.

\begin{lemma}\label{lemma_mult} For any multi-index ${\bf i} = (i_0,\dots, i_{m-1})$, with $i_0 < \dots < i_{m-1}$, consider the linear subspace $\Lambda_{\bf i} = \overline {{q}_{i_0} \dots q_{i_{m-1}}} \cong {\bf P}^{m-1}$ in ${\bf P}^{m+1}$. Then, there exists a quadric $Q_{\bf i}$ in $H^0 ({\bf P}^{m+1}, \mc{I}_{\Gamma_{m+1}} (2))$ containing $\Lambda_{\bf i}$.
\end{lemma}

\begin{proof} By symmetry, it suffices to prove the lemma for the multi-index ${\bf i} = (0,1,\dots,m-1)$. Let $\Lambda = \{ x_{n-1} = x_n = 0\} \cong {\bf P}^{n-2}$ be as in Section \ref{sec_cast}. Then, the projection of $\Lambda$ under $\pi$ is precisely $\Lambda_{\bf i}$. Define $Q_{\bf i} = \pi (Q_{n-m-1})$. By lemma \ref{lemma_proj}, $Q_{\bf i}$ is the projection of the scroll $\Sigma$. Since $Q_{n-m-1}$ contains $\Lambda$, it follows that $Q_{{\bf i}}$ contains $\Lambda_{{\bf i}}$.
\end{proof}

Finally, we have the important corollary, which should be thought of as an evidence towards Conjecture \ref{conj_2}.

\begin{cor}\label{cor_evidence} Suppose that for some (hence every) multi-index ${\bf i}$, the linear subspace $\Lambda_{\bf i}$ imposes independent conditions on quadrics containing $\Gamma_{m+1}$, i.e. the restriction map
$$\rho_{\bf i} : H^0 ({\bf P}^{m+1}, \mc{I}_{\Gamma_{m+1}}(2)) \rightarrow H^0 (\Lambda_{\bf i}, \mc{I}_{\Lambda_{\bf i} \cap \Gamma_{m+1}} (2))$$
has maximal rank.
 Then, Weak Conjecture holds for $\Gamma$.
\end{cor}
\begin{proof} We have: $$h^0 (\Lambda_{\bf i}, \mc{I}_{\Lambda_{\bf i} \cap {\Gamma_{m+1}}} (2)) = h^0 (\Lambda_{\bf i}, \mc{O}_{\Lambda_{\bf i}}(2)) - |\Lambda_{\bf i} \cap \Gamma_{m+1}| = \binom{m+1} 2 - m = \binom{m} 2.$$ 

By lemma \ref{lemma_mult}, $\rho_{\bf i}$ has a nontrivial kernel. Hence, $h^0 (\Lambda_{\bf i}, \mc{I}_{\Lambda_{\bf i} \cap {\Gamma_{m+1}}} (2)) \geq \binom{m} 2 + 1$. The claim follows from lemma \ref{lemma_tfae} and remark \ref{remark_inq}.
\end{proof}

\begin{remark} The corollary above gives a sufficient, but not a necessary condition for the Weak Conjecture to hold. For example, if $Bs(|\mc{I}_{\Gamma_{m+1}} (2)|)$ has a component $W$ of dimension $\geq 2$, then $W$ meets some $\Lambda_{\bf i}$ in a point that does not belong to $\Gamma_{m+1}$. In particular, $\rho_{\bf i}$ {\it fails} to be surjective in this case. 
\end{remark}

\bigskip

\section {Main Results} \label{sec_points2}

\subsection {The strategy}\

\bigskip

Finally, we are in position to present our results on Conjecture \ref{conj_2} for $m=3,4$. As before, let $\Gamma \subset {\bf P}^n$ be a set of $d \geq 2n+2m+1$ points in symmetric position ($m \leq n-2$), such that $h_\Gamma(2) = 2n+m$. We want to produce a curve $D$ of degree $\leq n+m-1$, containing $\Gamma$. Here is a naive approach, that will require a slight modification to work:
\bigskip
  
{\it Step 1.} Consider the projection $\pi_{m+1}: {\bf P}^n \dashrightarrow {\bf P}^{m+1}$ and the image $\Gamma_{m+1} \subset {\bf P}^{m+1} = \pi_{m+1} (\Gamma - \{p_0,\dots,p_{n-m-2}\})$. By using lemma \ref{lemma_mult}, we can show, that $\Gamma_{m+1}$ lies on sufficiently many quadrics in ${\bf P}^{m+1}$. From this information, we can deduce that $\Gamma_{m+1}$ lies on a curve $D_{m+1}$ of small degree, which is a component of $Bs(|\mc{I}_{\Gamma_{m+1}}(2)|)$. 
\bigskip

{\it Step 2.} Let $\Sigma \subset {\bf P}^n$ be the $m$-fold scroll constructed in sec. \ref{sec_cast} and let $\Sigma_{m+1} = \pi_{m+1} (\Sigma) \subset {\bf P}^{m+1}$.  By lemma \ref{lemma_proj}, the map $\pi|_\Sigma : \Sigma \dashrightarrow \Sigma_{m+1}$ is birational. Since $D_{m+1} \subset \Sigma_{m+1}$, we may consider the strict transform\footnote{Recall, that if $f: X \dashrightarrow Y$ is a birational map between varieties, and $Z \subset Y$ is a closed subset such that $f^{-1}$ is defined over a dense open subset $Z_0 \subset Z$, then {\it the strict transform} of $Z$ is, by definition, the closure of $f^{-1} (Z_0)$ in $X$. } $D \subset \Sigma$ of $D_{m+1}$ under $\pi_{m+1}$. 
\bigskip

Hopefully, the curve $D$ constructed in Step 2 is the one we are looking for. Unfortunately, we don't know the behaviour of $D$ at the points $\{p_0,\dots, p_{n-m-2}\}$, where $\pi_{m+1}|_\Sigma$ is not defined. A priori, $D$ could have arbitrary multiplicity $\mu$ at those points (including $\mu = 0$), and therefore we don't have control on the degree of $D$ or the fact that $D$ passes through all points of $\Gamma$. To remedy the situation, we introduce another step: 
\bigskip

{\it Step $1 \frac 1 2$.} Consider also the projection $\Gamma_{m+2} \subset {\bf P}^{m+2}$. By using the fact, that $\Gamma_{m+2}$ lies on sufficiently many quadrics, we can show that $\Gamma_{m+2}$ lies on a curve $D_{m+2}$ of small degree, which is a component of $Bs(|\mc{I}_{\Gamma_{m+2}}(2)|)$. Finally, we have the following ``lifting'' lemma, which exploits the symmetry of $\Gamma$ in a subtle way.

\begin{lemma}\label{lemma_lift} Let $\Gamma$ be as above. Assume, that $\Gamma_{m+2} \subset {\bf P}^{m+2}$ lies on a unique (reduced, irreducible) curve $D_{m+2}$ of degree $\leq 2m+2$. Let $D_{m+1} \subset {\bf P}^{m+1}$ be the projection of $D_{m+2}$ from $\pi_{m+2} (p_{n-m-2})$. Assume also, that $D_{m+1}$ is the unique curve in ${\bf P}^{m+1}$ of degree $\leq 2m+1$ that contains $\Gamma_{m+1}$.
Then, $\Gamma \subset {\bf P}^n$ lies on a curve $D$, such that $\pi_{m+1} (D) = D_{m+1}$. Moreover, $D$ is nonsingular at $\Gamma$ and $D$ has degree $(n-m-1)+\deg(D_{m+1})$.
\end{lemma}

\begin{proof} Notice, that $D_{m+2}$ is nonsingular at all points of $\Gamma_{m+2}$ (otherwise, by symmetry, $D_{m+2}$ is singular at all points of $\Gamma_{m+2}$; but, the genus of $D_{m+2}$ is $g(D_{m+2}) \leq m < |\Gamma_{m+2}|$). 

For any $k$, let $\Sigma_k = \pi_k(\Sigma) \subset {\bf P}^k$ be the projection of the $m$-fold rational normal scroll $\Sigma$, constructed in section \ref{sec_cast}. 
By cor. \ref{cor_birat}, the map $\pi_{k}|_{\Sigma} : \Sigma \dashrightarrow \Sigma_{k}$ is birational. 

Since 
\begin{align*}
|\Gamma_{m+2}| = d-(n-m-2) & \geq 2n+2m+1 - (n-m-2) \\
                           & \geq 4m+5 > 2 \deg (D_{m+2}),\end{align*}
$D_{m+2}$ is contained in $Bs (|\mc{I}_{\Gamma_{m+2}}(2)|)$. In particular, $D_{m+2} \subset \Sigma_{m+2}$. Similarly, $D_{m+1} \subset \Sigma_{m+1}$. 

Define $D \subset \Sigma$ to be the strict transform of $D_{m+2}$ with respect to the map $\pi_{m+2}|_{\Sigma} : \Sigma \dashrightarrow \Sigma_{m+2}$. A priori, $D$ is a curve in ${\bf P}^n$, containing $p_{n-m-2}, \dots, p_{d-1}$.
Since $D_{m+2}$ is nonsingular at $\Gamma_{m+2}$, $D$ is nonsingular at $p_{n-m-2},\dots,p_{d-1}$. 

Notice, that $D \subset \Sigma$ can be also defined as the strict transform of $D_{m+1}$ with respect to the map $\pi_{m+1}|_{\Sigma} : \Sigma \dashrightarrow \Sigma_{m+1}$. 

The key observation is that the later definition of $D$ does not depend on the order of points $p_0,\dots, p_{n-m-2}$ (In other words, let $K \subset L$ be the field extension as in Def. \ref{def_sym}. Then, $D$ is defined over an intermediate field $K \subset F \subset L$, such that $Gal(L/F)$ acts on $\{p_0,\dots, p_{n-m-2}\}$ via the full symmetry group on $n-m-1$ letters).

Since $D$ contains $p_{n-m-2}$ and is nonsingular at $p_{n-m-2}$, it follows by symmetry, that $D$ contains all of $p_0,\dots, p_{n-m-2}$ and $D$ is nonsingular at these points. 

In summary, $D$ is contains all of $\Gamma$ and is nonsingular at $\Gamma$.


Since $D_{m+1} = \pi_{m+1}(D)$, it follows, that $\deg(D) = (n-m-1) + \deg(D_{m+1})$.
\end{proof}

\subsection {The case $m=2$}\label{sec_m2}\

\bigskip

We reprove the well-known case $m=2$ of Conjecture \ref{conj_2} for points in {\it symmetric position}\footnote {In fact, the statement is true if $\Gamma$ is in {\it uniform position}, but our method does not handle the more general case.}. We use some auxiliary results from section \ref{sec_aux1} and \ref{sec_aux2} below.

\begin{prop}[Fano, Eisenbud-Harris] Let $\Gamma$ be a set of $d \geq 2n+5$ points in symmetric position in ${\bf P}^n$ ($n \geq 4$). Suppose that $h_\Gamma(2) = 2n+2$. Then, $\Gamma$ lies on an elliptic normal curve $D$.
\end{prop}

\begin{proof} We work in the setting of section \ref{sec_grob}. By cor. \ref{cor_mono}, the following degree 2 monomials marked with $''\bullet''$ and $''\blacksquare''$ appear in the initial ideal of $\Gamma$:
\begin{singlespace}
\begin{align*}
\begin{array}{lllllll}
\bullet & \bullet & \dots & \bullet & \bullet & \bullet & \circ \\
        & \bullet & \dots & \bullet & \bullet & \bullet & \circ \\
        &         & \ddots& \vdots  & \vdots  & \vdots  & \vdots\\
        &         &        & \bullet & \bullet & \bullet & \circ \\
        &         &        &         & \bullet & \blacksquare & \circ \\
        &         &        & &         & \circ & \circ \\
        &         &        & &         &         & \circ \\
\end{array}
\end{align*}
\end{singlespace}
\begin{center} {\bf The initial ideal of $\Gamma$ in degree 2 ($m=2,n\geq 4$)}
\end{center}

By cor. \ref{cor_codim}, there are no other initial monomials of degree 2 in $in(I_\Gamma)$, besides the ones on the diagram above.  

By lemma \ref{lemma_aux3} (in section \ref{sec_aux2} below), $\Gamma_4$ lies on an elliptic normal curve $D_4$. By the ``lifting'' lemma \ref{lemma_lift}, $\Gamma$ lies on an elliptic normal curve $D$.
\end{proof}

\subsection {The case $m=3$}\label{sec_m3}\

\bigskip

We prove the case $m=3$ of Conjecture \ref{conj_2} for points in symmetric position. We use some auxiliary results from section \ref{sec_aux1} and \ref{sec_aux2} below.

\begin{prop} Let $\Gamma$ be a set of $d \geq 2n+7$ points in symmetric position in ${\bf P}^n$ ($n \geq 5$). Suppose that $h_\Gamma(2) = 2n+3$. Then, $\Gamma$ lies on a curve $D$ of degree $\leq n+2$.
\end{prop}

\begin{proof}

We work in the setting of section \ref{sec_grob}. By cor. \ref{cor_mono}, the following degree 2 monomials marked with $''\bullet''$ and $''\blacksquare''$ appear in the initial ideal of $\Gamma$:
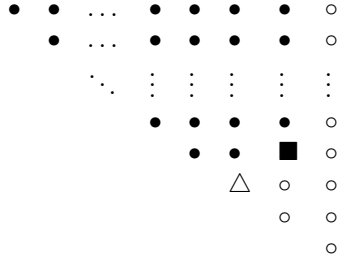
\begin{singlespace}
\begin{align*}
\begin{array}{llllllll}
\bullet & \bullet & \dots   &\bullet & \bullet & \bullet & \bullet & \circ \\
        & \bullet & \dots   &\bullet & \bullet & \bullet & \bullet & \circ \\
        &         & \ddots & \vdots  & \vdots  & \vdots  & \vdots  & \vdots\\
        &         &         &\bullet&\bullet  &\bullet  &\bullet  & \circ\\
        &         &         &        & \bullet& \bullet & \blacksquare & \circ \\
        &         &         &        &        & \triangle  & \circ & \circ \\
        &         &         &        &        &         & \circ & \circ \\
        &         &         &        &        &         &         & \circ \\
\end{array}
\end{align*}
\end{singlespace}
\begin{center} {\bf Part of the initial ideal of $\Gamma$ in degree 2 ($m=3,n\geq 5$)}
\end{center}

By cor. \ref{cor_codim}, only $\binom {3-1} 2 = 1$ monomial is missing from $in(I_\Gamma)_2$ in the diagram above.  

Consider the projection $\pi : {\bf P}^n \dashrightarrow {\bf P}^4$ and the image $\Gamma_4 \subset {\bf P}^4$ of $\Gamma - \{p_0,\dots,p_{n-5}\}$. By lemma \ref{lemma_mult} and lemma \ref{lemma_aux1} (in section \ref{sec_aux1} below), we have $h^0 ({\bf P}^4, \mc{I}_{\Gamma_4} (2)) \geq 4$. 

In particular, the missing monomial from $in(I_\Gamma)_2$ is the one marked with $''\triangle''$ in the diagram above.

By lemma \ref{lemma_aux4} (in section \ref{sec_aux2} below), $\Gamma_5$ lies on a curve $D_5$ of degree $\leq 7$. 

By the ``lifting'' lemma \ref{lemma_lift}, $\Gamma$ lies on a curve $D$ of degree $\leq n+2$.
\end{proof}

\subsection {The case $m=4 (n \geq 7)$}\label{sec_m4}\

\bigskip

We prove the case $m=4$, $n\geq 7$ of Conjecture \ref{conj_2} for points in symmetric position. We use some auxiliary results from section \ref{sec_aux1} and \ref{sec_aux2} below.

\begin{prop} Let $\Gamma$ be a set of $d \geq 2n+9$ points in symmetric position in ${\bf P}^n$ ($n \geq 7$). Suppose that $h_\Gamma(2) = 2n+4$. Then, $\Gamma$ lies on a curve $D$ of degree $\leq n+3$.
\end{prop}

\begin{remark} We are unable to handle the case $n = 6$, although the statement should still be conjecturally true. See the comment at the end of section \ref{sec_aux2} for an explanation.
\end{remark}

\begin{proof}

We work in the setting of section \ref{sec_grob}. By cor. \ref{cor_mono}, the following degree 2 monomials marked with $''\bullet''$ and $''\blacksquare''$ appear in the initial ideal of $\Gamma$:
\newpage

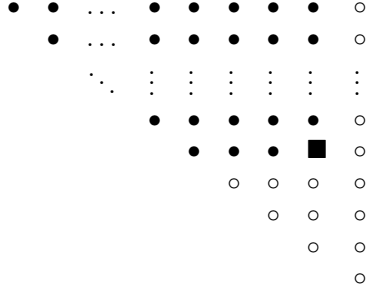
\begin{singlespace}
\begin{align*}
\begin{array}{lllllllll}
\bullet & \bullet & \dots  & \bullet   &\bullet & \bullet & \bullet & \bullet & \circ \\
        & \bullet & \dots & \bullet   &\bullet & \bullet & \bullet & \bullet & \circ \\
        &         &  \ddots & \vdots  & \vdots  & \vdots  & \vdots  & \vdots  & \vdots\\
        &         &         & \bullet &\bullet&\bullet  &\bullet  &\bullet  & \circ\\
        &         &         &         &\bullet&\bullet  &\bullet  &\blacksquare  & \circ\\
        &         &         &         &        & \circ  & \circ   & \circ & \circ \\
        &         &         &         &        &        & \circ   & \circ & \circ \\
        &         &         &         &        &        &         & \circ & \circ \\
        &         &         &         &        &        &         &         & \circ \\
\end{array}
\end{align*}
\end{singlespace}
\begin{center} {\bf Part of the initial ideal of $\Gamma$ in degree 2 (with $m=4,n\geq 7$)}
\end{center}

By cor. \ref{cor_codim}, only $\binom {4-1} 2 = 3$ monomials are missing from $in(I_\Gamma)_2$ in the diagram above.  

Consider the projection $\pi : {\bf P}^n \dashrightarrow {\bf P}^5$ and the image $\Gamma_5 \subset {\bf P}^5$ of $\Gamma - \{p_0,\dots,p_{n-6}\}$. By lemma \ref{lemma_mult} and lemma \ref{lemma_aux1} (in section \ref{sec_aux1} below), we have $h^0 ({\bf P}^5, \mc{I}_{\Gamma_5} (2)) \geq 6$. 

The possible missing monomials of degree 2 in the initial ideal $in(I_\Gamma)$ are depicted in the diagrams below (we set $n=7$ for simplicity).
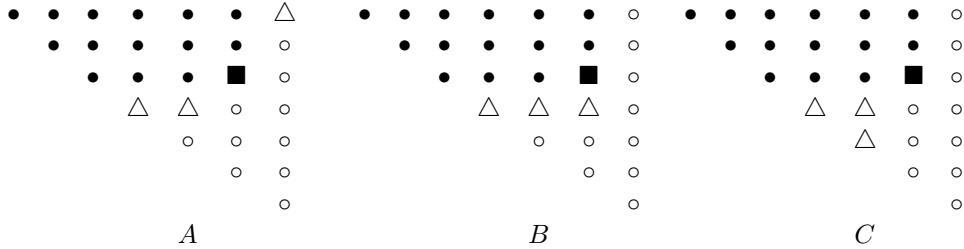
\begin{singlespace}
\begin{align*}
\begin{array}{ccccccc}
 \bullet& \bullet &\bullet &\bullet   & \bullet   & \bullet    & \triangle\\
        & \bullet &\bullet &\bullet    &\bullet   &\bullet    & \circ\\
        &         &\bullet &\bullet    &\bullet   &\blacksquare  & \circ\\
        &         &        & \triangle & \triangle& \circ    & \circ \\
        &         &        &           & \circ    & \circ     & \circ \\
        &         &        &           &          & \circ     & \circ \\
        &         &        &           &          &           & \circ \\
	&	  &	   &	       &  A       &	      &      \\
\end{array} \ \ \ \ 
\begin{array}{ccccccc}
 \bullet& \bullet &\bullet &\bullet   & \bullet   & \bullet    & \circ\\
        & \bullet &\bullet &\bullet    &\bullet   &\bullet    & \circ\\
        &         &\bullet &\bullet    &\bullet   &\blacksquare  & \circ\\
        &         &        & \triangle & \triangle& \triangle & \circ \\
        &         &        &           & \circ    & \circ     & \circ \\
        &         &        &           &          & \circ     & \circ \\
        &         &        &           &          &           & \circ \\
	&	  &	   &	       &  B	  &	      &      \\
\end{array}
\ \ \begin{array}{ccccccc}
 \bullet& \bullet &\bullet &\bullet   & \bullet   & \bullet    & \circ\\
        & \bullet &\bullet &\bullet    &\bullet   &\bullet    & \circ\\
        &         &\bullet &\bullet    &\bullet   &\blacksquare  & \circ\\
        &         &        & \triangle & \triangle& \circ & \circ \\
        &         &        &           &\triangle    & \circ     & \circ \\
        &         &        &           &          & \circ     & \circ \\
        &         &        &           &          &           & \circ \\
	&	  &	   &	       &   C	  &	      &      \\
\end{array}
\end{align*}
\end{singlespace}
\begin{center} {\bf Possible initial ideals for $I_\Gamma$ in degree 2 (with $m=4,n=7$)}
\end{center}

By lemma \ref{lemma_aux5} (in section \ref{sec_aux2} below), $\Gamma_6 \subset {\bf P}^6$ lies on a curve $D_6$ of degree $\leq 10$. 

By the ``lifting'' lemma \ref{lemma_lift}, $\Gamma$ lies on a curve $D$ of degree $\leq n+4$. (In particular, the diagram A for $in(I_\Gamma)_2$ is not possible, because $\Gamma$ satisfies the Weak conjecture, and so $x_0 x_n \notin in(I_\Gamma)_2$).

From the fact, that $h_D(2) = h_\Gamma(2) = 2n+4$, it follows, that $\deg D \leq n+3$ (this can be seen by taking a general hyperplane section). 

\end{proof}

\subsection {Auxiliary results I}\label{sec_aux1}\

\bigskip

The following elementary lemma on pencils of quadrics in projective space will be quite useful.

\begin{ellemma}\label{lemma_elt} Let $V \subset H^0 ({\bf P}^n, \mc{O}_{{\bf P}^n}(2))$ be a pencil of quadrics, containing a linear space $\Lambda = \overline{q_0\dots q_{n-3}} \cong {\bf P}^{n-3}$ in its base locus. Suppose that $V$ is generated by $Q_0,Q_1$. 
\begin{enumerate}
\item The locus of points $q \in {\bf P}^n$, for which there is a quadric in the pencil $V$ containing the span $\overline {\{q\}\cup \Lambda}$, is a subset of the determinantal locus
\begin{align}\label{eqn_loc} \text{rk} \left|
\begin{array}{cccc}
Q_0 (q_0, -) & Q_0 (q_1, -) & \dots & Q_0 (q_{n-3},-) \\
Q_1 (q_0, -) & Q_1 (q_1, -) & \dots & Q_1 (q_{n-3},-) \\
\end{array}
\right| \leq 1.
\end{align}
Here, we think of $Q_i(-,-)$ as a symmetric bilinear form on the underlying affine space ${\bf A}^{n+1}$.
\item Suppose that the determinantal locus in (\ref{eqn_loc}) defines the whole ${\bf P}^n$. Then, either $V$ contains a quadric singular along $\Lambda$ or else, every quadric in $V$ is singular along a fixed linear subspace ${\bf P}^{n-4} \subset \Lambda$.
\end{enumerate}
\end{ellemma}

\begin{proof} (a) Left to reader.

(b) Suppose that the determinantal locus in (\ref{eqn_loc}) defines the whole ${\bf P}^n$. Then, after some row and column operations, the matrix in (\ref{eqn_loc}) can be brought to one of the forms:
$$
\left|
\begin{array}{cccc}
0 & \dots & 0 & 0 \\
* & \dots & * & * \\
\end{array}
\right|
\ \ \ \ \text{or} \ \ \ \ \left|
\begin{array}{cccc}
0 & \dots & 0 & * \\
0 & \dots & 0 & * \\
\end{array}
\right|.
$$
The claim follows.
\end{proof}

{\bf Conventions.} In the next two lemmas, $\Gamma = \{ q_0, \dots, q_{d-1} \} \subset {\bf P}^n$ is a set of $d$ points in symmetric position. By {\it multi-index} ${\bf i} = (i_0,\dots,i_{k})$ we simply mean a $(k+1)$-tuple of distinct integers between $0$ and $d-1$. For any multi-index ${\bf i}$, we define $\Lambda_{\bf i} = \overline {q_{i_0} \dots q_{i_{k}}} \cong {\bf P}^k$. We work with homogeneous coordinates $x_0,\dots,x_n$ in ${\bf P}^n$, for which $q_0,\dots,q_n$ are the standard vertices. We consider the initial ideal of $I_\Gamma$ in the {\it lexicographic order}, as in section \ref{sec_borel}.

\begin{lemma}\label{lemma_aux1} Let $\Gamma = \{ q_0, \dots, q_{d-1} \} \subset {\bf P}^4$ be a set of $d \geq 16$ points in symmetric position. Suppose that for any multi-index ${\bf i} = (i_0,i_1,i_2)$, there is a quadric $Q_{\bf i}$ in $|\mc{I}_{\Gamma}(2)|$ containing $\Lambda_{\bf i} = \overline {q_{i_0} q_{i_1} q_{i_2}} \cong {\bf P}^2$. Then, $h^0 ({\bf P}^4, \mc{I}_\Gamma(2)) \geq 4$.
\end{lemma}

\begin{proof} We assume that $h^0 ({\bf P}^4, \mc{I}_\Gamma(2)) = 3$ and derive contradiction. 

For any ${\bf i} = (i_0,i_1,i_2)$, consider the restriction map:
$$\rho_{\bf i} : H^0 ({\bf P}^4, \mc{I}_\Gamma(2)) \rightarrow H^0 (\Lambda_{\bf i}, \mc{I}_{\Lambda_{\bf i} \cap \Gamma}(2)).$$
Since the map has a 1-dimensional kernel, it's image is a pencil of quadric in $\Lambda_{\bf i}$. Hence, the linear system $|\mc{I}_\Gamma(2)|$ has an extra base-point $r_{\bf i} \in \Lambda_{\bf i}$, besides $\Lambda_{\bf i} \cap \Gamma = \{q_{i_0},q_{i_1},q_{i_2}\}$. By symmetry, $r_{\bf i}$ does not lie on any side of the triangle with vertices $q_{i_0},q_{i_1},q_{i_2}$.

By prop. \ref{prop_borel}, the possible shapes for $in(I_\Gamma)_2$ are the following:
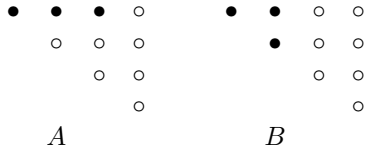
\begin{singlespace}
\begin{align*}
\begin{array}{cccc}
\bullet &\bullet &\bullet  &\circ   \\
        &\circ   &\circ    &\circ   \\
        &        &\circ    & \circ \\
        &        &         & \circ \\
        &  A     &         &       \\
\end{array} \ \ \ \ \ \
\begin{array}{cccc}
\bullet &\bullet &\circ    &\circ   \\
        &\bullet &\circ    &\circ   \\
        &        &\circ    & \circ \\
        &        &         & \circ \\
        &  B     &         &       \\
\end{array}
\end{align*}
\end{singlespace}
The shape {\it B} is ruled out, because there is no room for a quadric that contains $\Lambda_{012}$.

So, assume that $in(I_\Gamma)_2$ has shape A. Then, there is no quadric in ${\bf P}^4$ that both contains $\Gamma$ and is singular at $q_0$. By symmetry, there is no quadric that both contains $\Gamma$ and is singular at any point $q_i \in \Gamma$.

Now, the line $\overline {q_0 q_1}$ is not in the base-locus of the linear system $|\mc{I}_\Gamma(2)|$ (otherwise, by symmetry, any line through two points of $\Gamma$ is in the base-locus). Let $V$ be the pencil of quadrics that contain $\Gamma \cup \overline{q_0 q_1}$.

By lemma \ref{lemma_elt}(a), $\Gamma$ lies on the locus
\begin{align}\label{eqn_loc1} \text{rk} \left|
\begin{array}{ccc}
Q_0 (q_0, -) & Q_0 (q_1, -) \\
Q_1 (q_0, -) & Q_1 (q_1, -) \\
\end{array}
\right| \leq 1.
\end{align}
Now, the determinant of the matrix in eqn. (\ref{eqn_loc1}) must be identically 0, because there is no quadric that both contains $\Gamma$ and is singular at $\overline {q_0 q_1}$.

By lemma \ref{lemma_elt}(b), every quadric in the pencil $V$ is singular at a fixed point $s_{01} \in \overline {q_0 q_1}$. By symmetry, $s_{01}$ is distinct from $q_0$ and $q_1$. 

More generally, for any ${\bf j} = (j_0,j_1)$, we may consider the pencil of quadric that contain $\overline{q_{j_0} q_{j_1}}$. By symmetry, any quadric in that pencil is singular at a fixed point $s_{j_0 j_1} \in \overline {q_{j_0} q_{j_1}}$. 

For any multi-index ${\bf i} = (i_0,i_1,i_2)$, consider the plane $\Lambda_{i}$ and the points $q_{i_0},q_{i_1},q_{i_2},s_{i_0 i_1},$ $r_{i_0 i_1 i_2}$. It is easy to see, that
$$s_{i_0 i_1} = \overline {q_{i_0} q_{i_1}} \cap \overline {q_{i_2} r_{{i_0}{i_1}{i_2}}}$$
Since $r_{i_0,i_1,i_2}$ does not lie on any side of the triangle with vertices $q_{i_0}, q_{i_1}, q_{i_2}$, the points $s_{i_0 i_1}, s_{i_0 i_2}, s_{i_1 i_2}$ are not collinear.


Now, the quadric $Q_{012}$ is singular at the 3 noncollinear points $s_{01}, s_{02}, s_{12}$. Hence, $Q_{012}$ is singular along $\Lambda_{012}$. Contradiction!

This concludes the proof of the lemma.
\end{proof}

\begin{lemma}\label{lemma_aux2} Let $\Gamma = \{ q_0,\dots, q_{d-1} \} \subset {\bf P}^5$ be a set of $d \geq 20$ points in symmetric position. Suppose that for any multi-index ${\bf i} = (i_0,\dots,i_3)$, there is a quadric $Q_{\bf i}$ in $|\mc{I}_{\Gamma}(2)|$ containing $\Lambda_{\bf i} = \overline {q_{i_0}\dots q_{i_3}} \cong {\bf P}^3$. Then, $h^0 ({\bf P}^5, \mc{I}_\Gamma(2)) \geq 6$.
\end{lemma}

\begin{proof} We assume that $h^0 ({\bf P}^5, \mc{I}_\Gamma(2)) = 5$ and derive contradiction (the case $h^0 \leq 4$ is easier). By prop. \ref{prop_borel}, the possible shapes for $in(I_\Gamma)_2$ are the following:
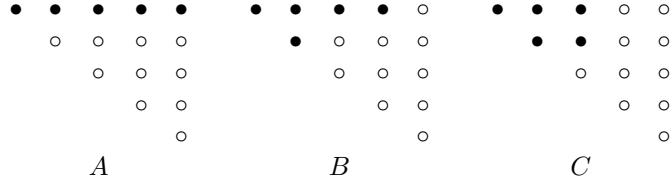
\begin{singlespace}
\begin{align*}
\begin{array}{ccccc}
 \bullet& \bullet &\bullet &\bullet  & \bullet   \\
        & \circ   &\circ   &\circ    &\circ   \\
        &         &\circ   &\circ    &\circ   \\
        &         &        &\circ    & \circ \\
        &         &        &         & \circ \\
        &         &  A     &         &       \\
\end{array} \ \ \ \ 
\begin{array}{ccccc}
 \bullet& \bullet &\bullet &\bullet  & \circ \\
        & \bullet   &\circ   &\circ    &\circ   \\
        &         &\circ   &\circ    &\circ   \\
        &         &        &\circ    & \circ \\
        &         &        &         & \circ \\
        &         &  B     &         &       \\
\end{array} \ \ \ \ 
\begin{array}{ccccc}
 \bullet& \bullet &\bullet &\circ    & \circ     \\
        & \bullet &\bullet &\circ    &\circ   \\
        &         &\circ   &\circ    &\circ   \\
        &         &        &\circ    & \circ \\
        &         &        &         & \circ \\
        &         &  C     &         &       \\
\end{array}
\end{align*}
\end{singlespace}
The shape {\it C} is ruled out, because there is no room for a quadric that contains $\Lambda_{0123}$.

From the description of $in(I_\Gamma)_2$, we see that there is no quadric that both contains $\Gamma$ and is singular along $\overline {q_0 q_1}$. By symmetry, there is no quadric that both contains $\Gamma$ and is singular along any line joining two points of $\Gamma$. 

Next, consider the number of conditions the 2-plane $\Lambda_{012}$ imposes on $|\mc{I}_\Gamma(2)|$. There are two conceivable cases. 

{\it Case I.} Suppose that $\Lambda_{012}$ imposes 3 independent conditions on $|\mc{I}_\Gamma(2)|$. In other words, there is a pencil of quadrics $V \subset H^0 ({\bf P}^5, \mc{I}_\Gamma (2))$, containing $\Lambda_{012}$. In particular, for any $q \in \Gamma$, there is a quadric in $V$, that contains the span $\overline { \{q\} \cup \Lambda_{012}}$. 

By lemma \ref{lemma_elt}(a), $\Gamma$ lies on the determinantal locus
\begin{align}\label{eqn_loc2} \text{rk} \left|
\begin{array}{ccc}
Q_0 (q_0, -) & Q_0 (q_1, -) & Q_0 (q_2,-) \\
Q_1 (q_0, -) & Q_1 (q_1, -) & Q_1 (q_2,-) \\
\end{array}
\right| \leq 1.
\end{align}
Now, all $2\times 2$ minors in (\ref{eqn_loc2}) must be identically 0, because there is no quadric that both contains $\Gamma$ and is singular along a line joining two points of $\Gamma$.

By lemma \ref{lemma_elt}(b), any quadric in the pencil $V$ must be singular along some fixed line $L_{012} \subset \Lambda_{012}$. By symmetry, $L_{012}$ does not pass through any of the points $q_0,q_1,q_2$.

More generally, for any multi-index ${\bf j} = (j_0,j_1,j_2)$, we can consider the pencil of quadrics in $|\mc{I}_\Gamma(2)|$ containing $\Lambda_{\bf j}$. By symmetry, any quadric in that pencil must be singular along a fixed line $L_{\bf j} \subset \Lambda_{\bf j}$. 

Let ${\bf i} = (i_0,i_1,i_2,i_3)$. Then, the quadric $Q_{\bf i}$ is singular along the 4 lines $L_{i_0 i_1 i_2}$, $L_{i_0 i_1 i_3}$, $L_{i_0 i_2 i_3}$, $L_{i_1 i_2 i_3}$. Since $Q_{\bf i}$ is irreducible (by linear generality of $\Gamma$), we conclude that the four lines span a 2-plane in $\Lambda_{\bf i}$.


Consider the projection $\pi : {\bf P}^5 \dashrightarrow {\bf P}^3$ from $L_{012}$ and let $\overline {\Gamma} \subset {\bf P}^3$ be the image of $\Gamma-\{q_0,q_1,q_2\}$. Then, $\overline{\Gamma}$ is in the base locus of the pencil of quadrics $\overline{V}$, which is the projection of the pencil $V$.

For any $j =3,\dots,d$, consider the quadric $Q_{012j}$ in $V$ and its projection $\overline {Q_{012j}}$ in $\overline V$. Since $Q_{012j}$ is singular along $L_{01j} \neq L_{012}$, the projection $\overline {Q_{012j}}$ is a singular quadric.

In summary, $\overline {V}$ is a pencil of quadrics in ${\bf P}^3$ with at least $d-3$ singular quadrics. It follows, that every quadric in $\overline {V}$ is singular at a fixed point of ${\bf P}^3$. So, $\overline{V}$ contains a reducible quadric. This contradicts the linear generality of $\Gamma$. 

{\it Case II.} Suppose that $\Lambda_{012}$ imposes only 2 independent conditions on $|\mc{I}_\Gamma(2)|$. In other words, the restriction map
$$H^0 ({\bf P}^5, \mc{I}_\Gamma(2)) \rightarrow H^0 (\Lambda_{012}, \mc{I}_{\\\Gamma \cap \Lambda_{012}} (2))$$
has rank 2. It follows, that the linear system $|\mc{I}_\Gamma(2)|$ has an extra base point $r_{012}$ in the plane $\Lambda_{012}$, besides $\Lambda_{012}\cap \Gamma = \{q_0,q_1,q_2\}$. By symmetry, $r_{012}$ is distinct from $q_0,q_1,q_2$.

By symmetry, for any multi-index ${\bf j} = (j_0,j_1,j_2)$, the linear system $|\mc{I}_\Gamma(2)|$ has an extra base-point $r_{\bf j}$ in the plane $\Lambda_{\bf j}$.

Let $W$ be the image of the restriction map
$$H^0 ({\bf P}^5, \mc{I}_\Gamma(2)) \rightarrow H^0 (\Lambda_{0123}, \mc{I}_{\\\Gamma \cap \Lambda_{0123}} (2)).$$
Since the map has 1-dimensional kernel, we have $\dim W = 4$. So, $W$ is a web of quadrics in ${\bf P}^3$ with at least 8 base-points $q_0,\dots,q_3, r_{012},\dots, r_{123}$. From what we know about the configuration of these 8 points, we easily conclude that such web $W$ cannot exist. Contradiction!

This completes the proof of the lemma.
\end{proof}

\begin{remark}\label{rmk_true} Assuming the hypotheses in the lemma above, we could naively expect the stronger conclusion $h^0 ({\bf P}^5, \mc{I}_{\Gamma} (2)) \geq 7$. Unfortunately, we don't know if this is true in general. In theory, this could fail if the linear system $|\mc{I}_\Gamma(2)|$ contains a ${\bf P}^2$ in its base locus (but, we think that such a possibility is unlikely). On the other hand, if it is true, then our proof of Conjecture \ref{conj_2} in the case $m =4$ could be simplified and extended to the case $n=6$.
\end{remark}

\subsection {Auxiliary results II}\label{sec_aux2}\

\bigskip

{\bf Conventions.} In the next three lemmas, $\Gamma = \{ q_0, \dots, q_{d-1} \} \subset {\bf P}^n$ is a set of $d$ points in symmetric position. We work with homogeneous coordinates $x_0,\dots,x_n$ in ${\bf P}^n$, for which $q_0,\dots,q_n$ are the standard vertices. We consider the initial ideal of $I_\Gamma$ in the {\it lexicographic order}, as in section \ref{sec_borel}.

\begin{lemma}\label{lemma_aux3} Let $\Gamma \subset {\bf P}^4$ be a set of $d \geq 13$ points in symmetric position. Suppose that $in(I_\Gamma)_2$ has the form:
\begin{singlespace}
\begin{align*}
\begin{array}{cccc}
\bullet &\bullet &\bullet  &\circ   \\
        &\bullet &\bullet  &\circ   \\
        &        &\circ    & \circ  \\
        &        &         & \circ  \\
\end{array}
\end{align*}
\end{singlespace}
In particular, $h^0 ({\bf P}^5, \mc{I}_{\Gamma} (2)) = 5$. Then, $\Gamma$ lies on an elliptic normal curve $D \subset {\bf P}^4$.
\end{lemma}

\begin{lemma}\label{lemma_aux4}  Let $\Gamma \subset {\bf P}^5$ be a set of $d \geq 17$ points in symmetric position. Suppose that $in(I_\Gamma)_2$ has the form:
\begin{singlespace}
\begin{align*}
\begin{array}{ccccc}
\bullet &\bullet &\bullet &\bullet & \circ \\
        &\bullet &\bullet &\bullet  &\circ   \\
        &        &\bullet &\circ    &\circ   \\
        &        &        &\circ    &\circ  \\
        &        &        &         &\circ  \\
\end{array}
\end{align*}
\end{singlespace}
In particular, $h^0 ({\bf P}^5, \mc{I}_{\Gamma} (2)) = 8$. Then, $\Gamma$ lies on a curve $D$ of degree $\leq 7$.
\end{lemma}

\begin{lemma}\label{lemma_aux5} Let $\Gamma \subset {\bf P}^6$ be a set of $d \geq 22$ points in symmetric position. Suppose that $in(I_\Gamma)_2$ has one of the forms:
\begin{singlespace}
\begin{align*}
\begin{array}{llllll}
\bullet &\bullet &\bullet    &\bullet   &\bullet    & \circ\\
        &\bullet &\bullet    &\bullet   &\bullet    & \circ\\
        &        &\bullet    &\bullet   & \circ     & \circ \\
        &        &           & \circ    & \circ     & \circ \\
        &        &           &          & \circ     & \circ \\
        &        &           &          &           & \circ \\
        &        &           &  A       &           &       \\
\end{array} \ \ \ 
\begin{array}{llllll}
\bullet &\bullet &\bullet    &\bullet   &\bullet    & \circ\\
        &\bullet &\bullet    &\bullet   &\bullet    & \circ\\
        &        &\bullet    &\bullet   &\bullet    & \circ \\
        &        &           & \circ    & \circ     & \circ \\
        &        &           &          & \circ     & \circ \\
        &        &           &          &           & \circ \\
        &        &           &  B       &           &       \\
\end{array} \ \ \ 
\begin{array}{llllll}
\bullet &\bullet &\bullet    &\bullet   &\bullet    & \circ\\
        &\bullet &\bullet    &\bullet   &\bullet    & \circ\\
        &        &\bullet    &\bullet   &\circ      & \circ \\
        &        &           &\bullet   &\circ      & \circ \\
        &        &           &          & \circ     & \circ \\
        &        &           &          &           & \circ \\
        &        &           &  C       &           &       \\
\end{array}
\end{align*}
\end{singlespace}
In particular, $h^0 ({\bf P}^6, \mc{I}_{\Gamma} (2)) \geq 11$. Then, $\Gamma$ lies on a curve $D$ of degree $\leq 10$.
\end{lemma}
\bigskip

Notice, that in all three lemmas, we have $x_0 x_n \notin in(I_\Gamma)_2$. Hence $\Gamma$ ``satisfies the Weak Conjecture'', i.e. for any point $p \in \Gamma$, the Zariski tangent space $T_p Bs(|\mc{I}_\Gamma(2)|)$ is 1-dimensional. Define $\Gamma'$ to be the zero-dimensional scheme supported on $\Gamma$, which contains the first infinitesimal neighborhood of $Bs(|\mc{I}_\Gamma(2)|)$ at each $p \in \Gamma$. In particular, $\deg \Gamma' = 2d$. 
\bigskip

{\it Proof of lemma \ref{lemma_aux3}.} In this case, $\deg \Gamma' \geq 26 > 2^4$. By Fulton's refinement of Bezout's theorem (\cite {Fult}), the base locus $Bs(|\mc{I}_\Gamma (2)|)$ contains a positive dimensional irreducible component $D$, passing through at least one point $p \in \Gamma$. Since the dimension of the Zariski tangent space to $D$ at $p$ is 1-dimensional, $D$ is a reduced curve. By symmetry, it follows that $D$ passes through all of $\Gamma$ (otherwise, $Bs(|\mc{I}_\Gamma(2)|)$ would have at least $d$ positive dimensional components, which is impossible). Since $D$ is a component of the intersection of 5 quadrics, we conclude that $D$ is an elliptic normal curve in ${\bf P}^4$. 
\bigskip

{\it Proof of lemma \ref{lemma_aux4}.} The argument is essentially the same as in the previous lemma. Here, we have $\deg \Gamma' \geq 34 > 2^5$, hence Bezout's theorem applies. As before, $\Gamma$ lies on a curve $D$, which is a component of $Bs (|\mc{I}_\Gamma(2)|)$. Since $D$ is a component of an intersection of 8 quadrics in ${\bf P}^5$, it follows, that $\deg D \leq 7$ (this can be seen by considering a general hyperplane section). 
\bigskip

{\it Proof of lemma \ref{lemma_aux5}.} This is the most technical lemma. Notice, that $\deg \Gamma' \geq 44 < 2^6$, so we cannot apply the same argument as above. 

{\it Step 1.} We claim, that for any point $q \in \Gamma$, there is a surface cone $C_q \subset {\bf P}^6$ with vertex $q$, containing $\Gamma'$, such that $\deg C_q \leq 10$.

It suffices to show this for just one point, say $q_0$. So, consider the projection $\Gamma_5 \subset {\bf P}^5$ of $\Gamma$ from the point $q_0$. From the shape of initial ideal of $\Gamma$, we see that $\Gamma_5$ still satisfies the Weak conjecture. Let $\Gamma_5'$ be the zero-dimensional scheme supported on $\Gamma_5$, containing the first infinitesimal neighborhood of $Bs(|\mc{I}_{\Gamma_5}(2)|)$ at each $\overline {q} \in \Gamma_5$. In particular, $\deg \Gamma'_5 = 2\deg \Gamma_5 \geq 42 > 2^5$. By the same argument as before, we conclude that $\Gamma'_5$ lies on a curve $Y$, which is a component of $Bs (|\mc{I}_{\Gamma_5}(2)|)$. From the shape of the initial ideal of $\Gamma$, $Y$ lies on at least 6 quadrics. It follows, that $\deg Y \leq 10$ (this can be seen by taking a general hyperplane section). Finally, let $C_{q_0}$ be the preimage of $Y$ in ${\bf P}^6$.

{\it Step 2.} We claim, that the base locus $Bs(|\mc{I}_\Gamma(2)|)$ contains a curve $D$, passing through a point of $\Gamma$. 

Consider the intersection of the cone $C_{q}$ with 2 general quadrics from $|\mc{I}_\Gamma(2)|$. Since $2^2 \cdot \deg C_{q} \leq 40$ and $\deg \Gamma' \geq 44$, we conclude that the intersection contains a positive dimensional component, passing through a point of $\Gamma$. Since the two quadrics are general, it follows that the base locus $Bs (|\mc{I}_\Gamma(2)|)$ contains a positive dimensional component $D$, passing through a point of $\Gamma$. Clearly, $D$ is a reduced curve.

{\it Step 3.} By symmetry, $D$ contains all of $\Gamma$. Since $D$ is a component of the intersection of $\geq 11$ quadrics in ${\bf P}^6$, it follows, that $\deg D \leq 10$ (again, this can be seen by considering a general hyperplane section). \hfill$\Box$


\newpage
\bibliographystyle{amsplain}

\begin{thebibliography}{10}


\bibitem {ACGH} {E. Arbarello, M. Cornalba, P. Griffiths, J. Harris}, {\it Geometry of Algebraic Curves, I,} Grundlehren der Mathematischen Wissenschaften, 267. Springer-Verlag (1985).

\bibitem {C} {G. Castelnuovo}, Ricerche di geometria sulle curve algebraiche, {\it Atti R. Accad. Sci Torino} {\bf 24} (1889), 196--223.

\bibitem {Ci} {C. Ciliberto}, Hilbert functions of finite sets of points and the genus of a curve in a projective space, {\it Lecture Notes in Math.} {\bf 1266} (1987), 24--73.

\bibitem {CLO} {D. Cox, J. Little, O'Shea}, {\it Ideals, Varieties and Algorithms. An introduction to computational algebraic geometry and commutative algebra}. Springer-Verlag, New York (1992).

\bibitem {Eis} {D. Eisenbud,} {\it Commutative algebra with a view toward algebraic geometry,} Graduate Texts in Math., 150. Springer-Verlag, New York (1995).

\bibitem {EGH} {D. Eisenbud, M. Green, J. Harris}, Higher Castelnuovo Theory, {\it Ast\'erisque} {\bf 218} (1993), 187--202.

\bibitem {EH} {D. Eisenbud, J. Harris}, {\it Curves in Projective Space}, Universit\'e de Montr\'eal (1982).

\bibitem {Fano} {G. Fano}, Sopra le curve di dato ordine e dei massimi generi in uno spazio qualunque, {\it Mem. Accad. Sci. Torino} {\bf 44} (1894), 335-382.

\bibitem{Fult} {W. Fulton}, {\it Intersection Theory}, Springer-Verlag, New York (1998).

\bibitem {G} {M. Green}, {Generic initial ideals}, {\it Six lectures on commutative algebra (Bellaterra, 1996)}. Progr. Math., {\bf 166}, Birh\"auser, Basel.

\bibitem {Gr} {L. Gruson, R. Lazarsfeld, C. Peskine}, {On a theorem of Castelnuovo and equations defining space curves}, {\it Invent. Math.} {\bf 72} (1983), no. 3, 491--506.

\bibitem {Halphen} {G. Halph\'en}, {\it M\'emoire sur la classification des courbes gauches,} in Oeuvres compl\`etes.

\bibitem {H} {J. Harris}, {\it A bound on the geometric genus of projective varieties}, Ph.D. Thesis, Harvard University (1978).

\end{thebibliography}

\end{document}